\documentclass[11pt]{article}
\usepackage[english]{babel}
\usepackage[latin1]{inputenc}
\usepackage{exscale}
\usepackage[centertags]{amsmath}
\usepackage{amsfonts}
\usepackage{amstext}		
\usepackage{amssymb}
\usepackage{amsthm}
\usepackage{newlfont}
\usepackage{enumerate}
\usepackage{graphics}
\usepackage{graphicx}
\usepackage{these}
\usepackage{color}
\title{The Isoperimetric Profile of a Compact Riemannian Manifold for Small Volumes\footnote{Work supported by a grant from INDAM}{}}
\author{Stefano Nardulli}
\begin{document}
      \maketitle
      \tableofcontents
      \addcontentsline{toc}{section}{\numberline{}Index}
      \addcontentsline{toc}{section}{\numberline{}Introduction}
      \newpage
\subsection{Isoperimetric profile, examples}
\begin{Def}
     Let $\mathcal{M}$ be a Riemannian manifold of dimension $n$ (possibly with infinite volume).\\
     Denote by $\tau_{\mathcal{M}}$ the set of  relatively compact open subsets of $\mathcal{M}$ where the  boundary is a submanifold of class        
     $C^{\infty}$.\\
     The function $I:[0,Vol(\mathcal{M})[\rightarrow [0,+\infty [$ such that
      $I(0)=0$
     \Fonct{I}{]0,Vol(\mathcal{M})[}{[0, +\infty [}{v}{Inf_{\left\lbrace \begin{array}{l}
                                                                                                       \Omega\in \tau_{\mathcal{M}}\\
                                                                                                       Vol(\Omega )=v 
                                                                                               \end{array} \right\rbrace }  \{ Vol_{n-1}(\partial \Omega )\} } 
is called the isoperimetric profile function (or shortly the isoperimetric profile) of the manifold $\mathcal{M}$.
\end{Def}
Here is a list of examples of Riemannian manifolds whose profiles are knowns.
\begin{enumerate}
          \item In constant sectional curvature spaces the solutions of the isoperimetric problem are geodesic     balls, this fact allow us, in the Euclidean case, to have an explicit formula for the isoperimetric profile and in the case of strictly positive or strictly negative curvature an implicit but satisfactory one.
         \item In dimension $2$, Frank Morgan, Michael Hutchings et Hugh Howards show in  
                  \cite{MHH} page 4899 (Main isoperimetric theorem) that on $\mathbb{R}^2$ endowed with a complete rotationally invariant metric whose sectional curvature is a strictly increasing fonction of the distance to the origin, the solutions of the isoperimetric problem are the following:
                  \begin{itemize}
                            \item a disk centered at the origin,
                            \item the complement of a disk centered at the origin,
                            \item a ring centered at the origin.
                  \end{itemize} 
                  This result generalises a previous result of Itai Benjamini and Jianguo Cao \cite{BC} page 361 (Main Theorem) that contains the case of the isoperimetric profile of a metric of paraboloids of revolution. Finally Manuel Ritor\'e in \cite{R} with different methods extends the list of surfaces for which the isoperimetric profile is known including planes of revolution with increasing sectional curvature etc. (see page 1094).
          \item Antonio Ros and Manuel Ritor\'e determine the isoperimetric profile of $\mathbb{R}P^3$   
                   in \cite{RR}, Laurent Hauswirth, Joaqu\'{i}n P\'erez, Pascal Romon and 
                   Antonio Ros in \cite{HPRR} do the same for certains flat tori of dimension $3$.
          \item Per Tomter \cite{Tom}  determines the profile only in a neighborhood of $0$, for the  
                   Heisenberg group of dimension $3$.  
          \item On certain surfaces of non trivial topology, one can construct metrics with somewhat prescribed isoperimetric profile, see \cite{GrPan}.  
\end{enumerate}
For a more detailed exposition of the state of the art we refer to the P.H.D. thesis of Vincent Bayle \cite{BaV}
and to the survey by Ros and Ritor\'e (on line on the web page of department of mathematics of the University of Granada). 
\subsection{Results}
The problem that we study in this article is the behaviour of the fonction $I_{\mathcal{M}}$ in a neighborhood of the origin. We compute the first nontrivial coefficient of its asymptotic expansion and open the way to the computation of other ones.

We show that the isoperimetric problem reduces, in a natural manner (in particular, compatible with the symmetry of the ambient manifold) to a variational problem in finite dimension. The geometric content of the following statement will become clear in subsequent paragraphs.

\begin{Res}
\label{1}
There exists a smooth function $f:\mathcal{M}\times\mathbb{R}\rightarrow\mathbb{R}$ and a map $\beta$ that to an arbitrary point $p$ of $\mathcal{M}$ and a real $r$ associates a domain $\beta(p,r)$, such that, for sufficiently small $r>0$, the solutions of the isoperimetric problem at volume $v=\omega_{n}r^{n}$ are exactly the images via $\beta$ of the minima of the function $p\mapsto f(p,r)$.  $f$ and $\beta$ are invariant (resp. equivariant) under the group of isometries of $\mathcal{M}$.
\end{Res}

As a consequence of theorem \ref{1}, in \cite{Tom}, Per Tomter determines the isoperimetric profile of the Heisenberg group of dimension $3$ endowed with a left invariant Riemannian metric of maximal symmetry, for small volumes. Tomter uses an unpublished work of Bruce Kleiner, dating back to 1985, that we can state as follows.

\begin{Thm}[Kleiner] \label{KleinerThm} Let $\mathcal{M}$ be a Riemannian manifold, $G$ a group of isometries that acts transitively on $\mathcal{M}$, $K\leq G$ the stabilizer of a point.\\
Then for small $v$, there exists a solution of the isoperimetric profile in volume $v$ that is $K$-invariant.  
\end{Thm}
We first learned of this result in may 2005, because there is no written trace and no other reference than  \cite{Tom}.

One of the goals of this paper is in particular to recover the proof of theorem \ref{KleinerThm}.
With some generic hypotheses on the Riemannian metric, we can make more precise the statement of theorem $1$ and we can show  the differentiability properties of the isoperimetric profile in a neighborhood of  $0$, for details see section 3.

\bigskip

The basic idea of the proof is to apply the implicit function theorem to the map that to a small hypersurface, seen like a small perturbation of a small geodesic sphere, associates its mean curvature. The linearization of this map is not invertible. Hence we modify this map, introducing a new class of hypersurfaces satisfying a weaker condition than the constancy of mean curvature. We call this new kind of hypersurfaces {\em pseudo-balls}. On the other hand it will be nedeed to show that a solution of  the isoperimetric problem is a perturbation of a small geodesic sphere.

The proof of theorem \ref{1} has two essentials steps,
\begin{enumerate}
          \item the construction of pseudo balls,
          \item the proof that the solutions of the isoperimetric profile are pseudo balls.
\end{enumerate}  
The second step is the subject of a separate article \cite{Nar}.
\subsection{Pseudo-balls}
\begin{Def}
          We call \textbf{pseudo-ball} a hypersurface $\mathcal{N}$ embedded in $\mathcal{M}$ such that there exists a point $p\in\mathcal{M}$ 
          and a function $u\in C^{2,\alpha}(T^1_p \mathcal{M}\backsimeq\mathbb{S}^{n-1},\mathbb{R})$, such that $\mathcal{N}$ is the graph of
$u$ in normal polar coordinates centered at $p$, i.e. $\mathcal{N}=\left\lbrace exp_p (u(\theta )\theta ),\:\theta\in T^1_p \mathcal{M} \right\rbrace$ and 
$$Q(H(u))=const.\in\mathbb{R},$$
where $H$ is the mean curvature operator, $Q=id-P$, $P$ is the orthogonal projector of 
$L^2 (T^1 _p \mathcal{M})$ on the first eigenspace of the Laplacian on the unit sphere $\mathbb{S}^{n-1}$ of $\mathbb{R}^{n}$, and $T^1 _p \mathcal{M}$ is the fiber over $p$ of the unit tangent bundle over the Riemannian manifold $\mathcal{M}$. 
\end{Def}
To state a uniqueness theorem for pseudo-balls we need the notion of \textit{center of mass}.
\begin{Def}
Let $(\Omega ,\mu )$ be a probability space and $f:\Omega\rightarrow\mathcal{M}$ a measurable function.
we consider the following function $\mathcal{E}:\mathcal{M}\rightarrow [0,+\infty[$:
$$\mathcal{E}(x):=\frac{1}{2}\int_{\Omega } d^2 (x,f(y))d\mu (y).$$
We call {\em center of mass} of $f$ with respect to the measure $\mu$ the unique minimum of 
$\mathcal{E}$ on $\mathcal{M}$, provided that it exists.
\end{Def}
\begin{Prop}\label{p1}
We assume that the sectional curvature of $\mathcal{M}$ satisfies $\mathcal{K}_{\mathcal{M}}\leq\delta$. We assume also that $f(\Omega)$ has a diameter less than one half of the injectivity radius and than $\frac{\pi}{2\sqrt{\delta}}$.\\
Then $\mathcal{E}$ has a unique minimum $c$, i.e. the unique nearest point to $f(\Omega)$ such that  
$\int_{\Omega}exp_c ^{-1}f(y)d\mu (y)=0$. 
\end{Prop}
In particular, we can speak about the center of mass of a hypersurface of small diameter (we apply proposition \ref{p1} to the $(n-1)$-dimensional measure of the boundary).
\begin{Def}
          Let $\mathcal{F}^{2,\alpha}$ be the fiber bundle on $\mathcal{M}$ with fiber at $p$ equal to  
          $C^{2,\alpha}(T^1_p \mathcal{M},\mathbb{R})$
\end{Def}  
\begin{Res}\label{TT1}
There exists a $C^{\infty}$ map, $\beta:\mathcal{M}\times\mathbb{R}\rightarrow\mathcal{F}^{2,\alpha}$ such that for all $p\in\mathcal{M}$, and all sufficiently small $r>0$, the hypersurface $exp_p (\beta(p,r)(\theta )\theta) $ is the unique pseudo-ball whose center of mass is $p$ enclosing a volume $\omega_n r^n$ where $\omega_n :=Vol_{(\mathbb{R}^n,can)}(\mathbb{B}^n )$ is the Euclidean volume of the unit ball $\mathbb{B}^n$ of $ \mathbb{R}^n$.
\end{Res}
\textbf{Remark:} If $g$ is an isometry of $\mathcal{M}$, $g$ sends pseudo-balls to pseudo-balls
and $g\circ\beta =\beta\circ g$ ($g$ acts only on the first factor $\mathcal{M}$).

\textbf{Examples:} 
\begin{enumerate}
          \item In constant curvature, pseudo-balls are geodesics spheres.
          \item In the $3$ dimensional Heisenberg group, pseudo-balls have been calculated by Per Tomter \cite{Tom}. 
\end{enumerate}
\subsection{Geometric measure theory}
We use the following result, from \cite{Nar}.
\begin{Res}\label{T4}
            Let $\mathcal{M}^n$ be a compact Riemannian manifold, $g_j$ a sequence of Riemannian metrics of  class $C^{\infty}$ that converges $C^4$ to a fixed metric $g_{\infty}$.
Let $B$ be a domain in $\mathcal{M}$ with smooth boundary $\partial B$, $T_j$ a solution of the isoperimetric problem in $(\mathcal{M}^n , g_j )$ such that  
$$(*): \mathbf{M}_{g_{\infty}} (B-T_j)\rightarrow 0.$$ 
Then $\partial T_j$ is the graph of a function $u_j$ in normal exponential coordinates to $\partial B$.
Furthermore for all $\alpha\in ]0,1[$, $u_j \in C^{2,\alpha}(\partial B)$, and if we assume that $\partial B$ has constant mean curvature for $g_{\infty}$, then $||u_j||_{C^{2,\alpha}(\partial B)}\rightarrow 0$.
\end{Res}
In section $3$, we deduce the following theorem, that completes the proof of theorem $1$.
\begin{Res}\label{T3}
 Let $T$ be a  current solution of the isoperimetric problem, of sufficiently small volume. Then $T$ is a pseudo-ball.
\end{Res}
\subsection{Expansion of the isoperimetric profile near $0$}
The asymptotic expansion of the volume of pseudo-balls and the volume of their boundary can be computed with theorem $1$, this yields an expansion for the profile.\\
The paper \cite{Druet} of O. Druet provides an alternative approach to theorem $7$. But our method gives more information. For instance, when scalar curvature is constant, the next term in the expansion of the profile can be obtained.   
\subsection{Plan of the article}
\begin{enumerate}
           \item Section $1$ describes the construction of pseudo balls via the implicit function theorem in an infinite dimensional context.
           \item Section $2$ describes why and in what sense approximate solutions, of the isoperimetric problem, in the case of small volumes, are close to Euclidean balls. 
           \item In section $3$ the results of preceding sections and those of \cite{Nar} are applied to obtain  the first two non zero coefficients of the asymptotic expansion of the isoperimetric profile.
           \end{enumerate}
\subsection{Aknowledgements} I wish to thank Renata Grimaldi of the University of Palermo for the fruitful discussions that we had during my P.H.D. studies, about the subject of my thesis. I thank also my P.H.D. advisor Pierre Pansu. I thank also Guy David, Séverine Rigot, Giandomenico Orlandi for their suggestions about the theorem of Jean Taylor. I am grateful to Istituto Nazionale di Alta Matematica "Francesco Severi" for financial support. 

      \newpage
      \section{Construction of pseudo-balls}
\subsection{The normal coordinates context} 
We place ourselves in $T_p \mathcal{M}$, the tangent space of $\mathcal{M}$ at $p$, endowed with the Riemannian metric $exp_p ^* (g)$. This is a smooth Riemannian metric in a neighborhood of the origin, $||x||<\varepsilon$. Let $u$ be a function on $T^1 _p \mathcal{M}\cong\mathbb{S}^{n-1}$ the unit sphere of $T_p \mathcal{M}$, such that $||u||_{L^{\infty}}<\varepsilon$, we are interested in the hypersurface 
$$\mathcal{N}_u =\{exp_p (u(\theta )\theta )|\theta\in T^1 _p \mathcal{M} \}.$$\\ 
Such hypersurfaces will be called \textit{normal graphs}.
Denote by $\theta$ the radial unit vector field on $T_p (\mathcal{M})-\{ 0\}$.
\subsection{The mean curvature}
Let $\mathcal{N}$ be an arbitrary hypersurface. It's mean curvature is
        \begin{equation}\label{h-intr} 
                        H_\nu ^{\mathcal{N}}(x)=\underset{1}{\overset{n-1}\sum }_{i}II^\mathcal{N}_x (e_i ,e_i )=
                                                -\underset{1}{\overset{n-1}\sum }_{i} <\nabla_{e_i }\nu , e_i >_g (x)
        \end{equation}
        where $ (e_1 ,\ldots ,e_{n-1} )$ is an  orthonormal basis of $T_x \mathcal{N}$, $\nu$ a normal unit vector field to $\mathcal{N}$ and  $II^\mathcal{N}_x$ the second fundamental form of 
$\mathcal{N}$ in $x$.\\
We will use in the sequel an extension of $\nu$ to the whole $T_p (\mathcal{M})-\{ 0\}$ into a vector field that is independent from the distance to the origin,
i.e. such that $[\theta, \nu]=0$ where $[,]$ is the Lie bracket of two vector fields.\\
We set
$$ \nu =a+b\theta ,$$
where the  vector field $a$ is tangent to geodesic spheres $S(p,r)$ centered at $p$ and satisfies obviously $ [a,\theta ]=0$.
In this case, after calculations analogous to what can be found in paragraph 3.3 of \cite{Nar} one gets:
\begin{equation}\label{cmpp}
      H_\nu (r,\theta)=-div_{S(p,r)}(a)+<\nabla_a a,a>-bII_{\theta}^r (a,a)+bH_{\theta}^r+b\nabla_a b
\end{equation} 
where $div_{S(p,r)}(a)$ is the divergence of the vector field $a$ restricted to the geodesic sphere of radius $r$ centered at $p$, $II_{\theta}^r$ is the second fundamental form of the sphere of radius $r$ centered at $p$ in the outward radial direction and $H_{\theta}^r$ is the trace of $II_{\theta}^r$.\\      
The next task is to reformulate (\ref{cmpp}) in terms of objects that live on the sphere $T^1_p \mathcal{M}$. 
We denote $\tilde{a}(r):=(di_r )^{-1} (a)$ if $a\in T_{exp_p(r\theta)}\mathcal{M}$, $$i_r:\left\{
             \begin{array}{lll}
                   T^1_p \mathcal{M} & \rightarrow & \mathcal{M}\\
                   \theta & \mapsto & exp_p (r\theta),\\ 
             \end{array}\right .
$$
we consider the pull-back $g_r$ of the metric $g$ with respect to $i_r$ (i.e. $g_r =i_r ^* (g)$). 
\begin{equation}\label{cmpp1}
      H_\nu =-div_{(\mathbb{S}^{n-1},g_r)}(\tilde{a})+<\nabla_{\tilde{a}} r\tilde{a},\tilde{a}>_{g_r}-br^2 II_{\theta}^r (\tilde{a},\tilde{a})+bH_{\theta}^r+b\nabla_{r\tilde{a}}b
\end{equation} 
We introduce the family of metrics $\tilde{g}_r :=\frac{1}{r^2}g_r$ that have the advantage to have no singularities in $r=0$, to highlight the nature of the singularity in $0$ of (\ref{cmpp}), thus (\ref{cmpp1}) becomes
 \begin{equation}\label{cmpp2}
      H_\nu =-div_{(\mathbb{S}^{n-1},\tilde{g}_r)}(\tilde{a})+r^2 <\nabla_{\tilde{a}} r\tilde{a},\tilde{a}>_{\tilde{g}_r}-br^2 II_{\theta}^r (\tilde{a},\tilde{a})+bH_{\theta}^r+b\nabla_{r\tilde{a}}b
\end{equation} 
\subsection{The mean curvature of a normal graph}
In this paragraph we interpret formulae (\ref{cmpp1}) and (\ref{cmpp2}) in the particular case of hypersurfaces $\mathcal{N}_u$. Let
$$W_u :=\sqrt{1+\|\overrightarrow{\nabla}_{i_u ^* (g)}u \| _{i_u ^* (g)}^2 }.$$ 
Then 
$$b= \left\{ 
         \begin{array}{llll}
                      \frac{1}{W_u } & <\nu , \theta >\geq 0 & \nu & outward\\
                     -\frac{1}{W_u } & <\nu , \theta >\leq 0 & \nu & inward\\
         \end{array}\right. 
         $$
\begin{equation}
           a=-bu\mathcal{I}_{\theta , u(\theta )\theta }\left( \overrightarrow{\nabla}_{g_u }u\right) 
\end{equation} 
where $\mathcal{I}_{\theta , u(\theta )\theta }$ is the identification of $\mathbb{R}^n\cong T_p \mathcal{M}$ of $T_{\theta}T_p \mathcal{M}$ with $T_{u(\theta )\theta }T_p \mathcal{M}$ that is induced by  the chosen system of normal coordinates centered in $p$ (notation inspired by \cite{Chav}).
In the sequel we always make the choice of $b=-\frac{1}{W_u }$, hence
\begin{equation}
          a=\frac{u}{W_u }\mathcal{I}_{\theta , u(\theta )\theta }\left( \overrightarrow{\nabla}_{g_u }u\right). 
\end{equation}         
For $r$ sufficiently small, formula (\ref{cmpp1}) becomes
\begin{eqnarray}\label{h-norm1}  
                    H_{\nu_{inw}} ^\mathcal{N}(u,\theta)&= &
                                   -div_{(\mathbb{S}^{n-1},g_u )}(\frac{\overrightarrow{\nabla}_{g_u} u}{W_{u}})
                               -\frac{1}{W_{u} ^2}<\nabla_{\overrightarrow{\nabla}_{g_u} u}(\frac{u\overrightarrow{\nabla}_{g_u} 
                                    u}{W_{u}}),\overrightarrow{\nabla}_{g_u} u>_{g_u}\\ \nonumber
                                   &+&\frac{u^2 }{W_{u} ^3}II_\theta ^{u} (\overrightarrow{\nabla}_{g_u} u,
                                     \overrightarrow{\nabla}_{g_u} u)
                                   -\frac{1}{W_{u}}H_\theta ^{u} (u,\theta)\\ \nonumber
                                   &+&\frac{1}{W_{u}}<\overrightarrow{\nabla}_{g_u} (\frac{1}{W_{u}}),
                                    u\frac{\overrightarrow{\nabla}_{g_u} u}{W_{u}})>_{g_u}.
\end{eqnarray}
Formula (\ref{cmpp2}) becomes
\begin{eqnarray}\label{h-norm2}  
                    H_{\nu_{inw}} ^\mathcal{N}(u,\theta)&= &
                                   -div_{(\mathbb{S}^{n-1},\tilde{g}_u )}(\frac{\overrightarrow{\nabla}_{\tilde{g}_u} u}{u^2 W_{u}})
                                -\frac{1}{u^2 W_{u} ^2}<\nabla_{\overrightarrow{\nabla}_{\tilde{g}_u} u}(\frac{\overrightarrow{\nabla}_{\tilde{g}_u} 
                                    u}{uW_{u}}),\overrightarrow{\nabla}_{\tilde{g}_u} u>_{\tilde{g}_u}\\ \nonumber
                                   &+&\frac{1}{u^2 W_{u} ^3}II_\theta ^{u} (\overrightarrow{\nabla}_{\tilde{g}_u} u,
                                     \overrightarrow{\nabla}_{\tilde{g}_u} u)\\ \nonumber
                                   &-&\frac{1}{W_{u}}H_\theta ^{u} (u,\theta)
                                   +\frac{1}{W_{u}}<\overrightarrow{\nabla}_{\tilde{g}_u} (\frac{1}{W_{u}}),
                                    \frac{\overrightarrow{\nabla}_{\tilde{g}_u} u}{uW_{u}})>_{\tilde{g}_u}.
\end{eqnarray}
\subsection{Linearization of the modified mean curvature}
        In this paragraph we denote by $T^1_p \mathcal{M}=\mathbb{S}^{n-1}T_p \mathcal{M}$ the unit sphere of $T_p \mathcal{M}$.
         \begin{Def}
          Let $\mathcal{F}^{k,\alpha}$ be the fiber bundle on $\mathcal{M}$ where the fiber over $p$ is 
          the space of $C^{k,\alpha}$ functions on the unit tangent sphere $T^1_p \mathcal{M}$ and $\Gamma (\mathcal{F}^{k,\alpha})$ the topological space of $C^{\infty}$ sections of $\mathcal{F}^{k,\alpha}$.
         \end{Def} 
      In other words, if $y\in\Gamma (\mathcal{F}^{k,\alpha})$, $p\in\mathcal{M}$, then $x=y(p)$ is a function on $T^1_p \mathcal{M}$. Let $y_0$ denote the zero section.\\
         In the rest of this section we are interested in the local behaviour of certain functions on this fiber bundle in view of the application of the implicit function theorem to the solution of equations in a neighborhood of the zero section $y_0$. In order to do this, it is convenient to denote $y=(p,x)$ and to identify a neighborhood in the fiber bundle $\mathcal{F}^{k,\alpha}$ with the trivialization 
$\mathcal{U}\times C^{2,\alpha}(\mathbb{S}^{n-1},\mathbb{R})$ ($\mathcal{U}$ open set of  $\mathcal{M}$ ) with the aid of an atlas of the differentiable structure of $\mathcal{M}$.\\   
         We define, now, the domains and codomains of the functionals of modified curvature.
       \begin{Def}We let
\Fonct{\Psi}{\mathbb{R}\times\Gamma(\mathcal{F}^{2,\alpha})}{\Gamma(\mathcal{F}^{0,\alpha})}{(r,y)}{r\left( H(p,r(1+x))-\frac{n-1}{r}\right).}     
In other words, if $p\in\mathcal{M}$, $r\in\mathbb{R}$ and $x$ is a $C^{2,\alpha}$ function on $T^1_p \mathcal{M}$, then $\Psi (p,r,x)$ is the $C^{0,\alpha}$ function on $T^1_p \mathcal{M}$ defined by              
                   $$\Psi (p,r,x):=r\left( H(p,r(1+x))-\frac{n-1}{r}\right).$$
       \end{Def}
        \begin{Prop}\label{ord1}$\Psi$ is $C^{\infty}$ on the open subset where $||x||_{\infty}<1$. 
(here $||\cdot||_{\infty}$ is the $L^{\infty}$-norm) 
        \end{Prop}
\begin{Dem}
           We put $u=r(1+x)$ in the formula (\ref{h-norm2}) and we get
\begin{enumerate}
           \item $W_{r(1+x)}=\sqrt{1+\frac{1}{r^2 (1+x)^2}||\overrightarrow{\nabla}_{\tilde{g}_{r(1+x)}}r(1+x)||_{\tilde{g}_{r(1+x)}}^2}=h_0(p,r,x)$,\\
                    $h_0(p,0,x)=\sqrt{1+\frac{1}{(1+x)^2 }||\overrightarrow{\nabla}_{\tilde{g}_0 }x||_{\tilde{g}_0 }^2}$
           \item $-div_{(\mathbb{S}^{n-1},\tilde{g}_{r(1+x)})}(\frac{\overrightarrow{\nabla}_{\tilde{g}_{r(1+x)}}r(1+x) }{r^2 (1+x)^2 W_{r(1+x)}}) =\frac{1}{r}h_1(p,r,x)$,\\ $h_1(p,0,x)=-div_{(\mathbb{S}^{n-1},can=g_0 )}(\frac{\overrightarrow{\nabla}_{\tilde{g}_{0}}x}{(1+x)^2 h_0 (p,0,x)})$ 
           \item $\frac{1}{r^2 (1+x)^2 W_{r(1+x)} ^2}<\nabla_{\overrightarrow{\nabla}_{\tilde{g}_{r(1+x)}} rx}(\frac{\overrightarrow{\nabla}_{\tilde{g}_{r(1+x)}} 
                                    rx}{r(1+x)W_{r(1+x)}}),\overrightarrow{\nabla}_{\tilde{g}_{r(1+x)}} rx>_{\tilde{g}_{r(1+x)}}=h_2(p,r,x)$,\\ $h_2 (p,0,x)=\frac{1}{(1+x)^2 h_0 (0,x)^2 }<\nabla_{\overrightarrow{\nabla}_{\tilde{g}_ 0 } x}(\frac{\overrightarrow{\nabla}_{\tilde{g}_0 } 
                                    x}{(1+x)h_0 (p,0,x)}),\overrightarrow{\nabla}_{\tilde{g}_0 } x>_{\tilde{g}_0 }$
            \item $\frac{1}{r^2 (1+x)^2 W_{r(1+x)} ^3}II_\theta ^{r(1+x)} (\overrightarrow{\nabla}_{\tilde{g}_{r(1+x)}} rx,
                                     \overrightarrow{\nabla}_{\tilde{g}_{r(1+x)}} rx)=\frac{1}{r}h_3(p,r,x)$,\\
            $h_3(p,0,x)=\frac{1}{(1+x)^3 h_0 (p,0,x)^3}(\overrightarrow{\nabla}_{\tilde{g}_0 }x,\overrightarrow{\nabla}_{\tilde{g}_0 }x)_{\mathbb{R}^{n-1}}$
            \item $-\frac{1}{W_{r(1+x)}}H_\theta ^{r(1+x)} (r(1+x),\theta)=\frac{1}{r}h_4(p,r,x)$\\
                     $h_4(p,0,x)=\frac{n-1}{(1+x)h_0 (p,0,x)}$
            \item $\frac{1}{W_{r(1+x)}}<\overrightarrow{\nabla}_{\tilde{g}_{r(1+x)}} (\frac{1}{W_{r(1+x)}}),\frac{\overrightarrow{\nabla}_{\tilde{g}_{r(1+x)}} rx}{r(1+x)W_{r(1+x)}})>_{\tilde{g}_{r(1+x)}}=h_5 (p,r,x)$,\\
                     $h_5 (p,0,x)=\frac{1}{h_0 (p,0,x)^2}<\overrightarrow{\nabla}_{\tilde{g}_0 } (\frac{1}{h_0 (p,0,x)}),\frac{\overrightarrow{\nabla}_{\tilde{g}_0 } x}{(1+x)})>_{\tilde{g}_0}$. 
\end{enumerate}
The functions $h_0$, $h_1$, $h_2$, $h_3$, $h_4$, $h_5$ are in $C^{\infty}(\mathbb{R}\times\Gamma(\mathcal{F}^{2,\alpha}),\Gamma(\mathcal{F}^{0,\alpha}))$ provided $||x||_{\infty}<1$.
From the following formulae
\begin{equation}
           rH(p,r(1+x))=h_1 (r,x)+rh_2 (r,x)+h_3 (r,x)+h_4 (r,x)+rh_5 (r,x),
\end{equation}
\begin{equation}
            \Psi (p,r,x)=rH(p,r(1+x))-(n-1),
\end{equation}  
we get that $\Psi$ is $C^{\infty}$.
\end{Dem}
        \\
        We can compare proposition \ref{ord1} with the calculation of the papers \cite[page 383]{Ye} and  \cite{Fethi}.\\ 
        \begin{Lemme}\label{lc11}  Let $P$ be the orthogonal projector of 
$L^2 (T^1 _p \mathcal{M})$ on the first eigenspace of the Laplacian on the unit sphere $\mathbb{S}^{n-1}$ and let $Q$ be $Id-P$. Denote by 
\Fonct{L}{C^{k+2,\alpha}(\mathbb{S}^{n-1})}{C^{k,\alpha}(\mathbb{S}^{n-1} )}{v}{-\triangle_{ _{\mathbb{S}^{n-1}}}v-(n-1)v} 
Then   $\left[\frac{\partial}{\partial x}Q(\Psi(p,r,x))\right] _{r=0, x=0}=L$.\\     
Furthermore, let $C_1 ^{l,\alpha}(\mathbb{S}^{n-1} ):=Ker(L)^{\perp}\cap C^{l,\alpha}(\mathbb{S}^{n-1} )$ where $Ker(L)^{\perp}$ is taken in the $L^2$ sense. Then 
$L:C_1^{k+2,\alpha}(\mathbb{S}^{n-1} )\mapsto C_1^{k,\alpha}(\mathbb{S}^{n-1} )$ is an isomorphism.         
        \end{Lemme} 
        \begin{Dem}
                  The following straightforward calculation shows that $$L(v)=-\triangle_{ _{\mathbb{S}^{n-1}}}v-(n-1)v.$$
                  With the notations of proposition \ref{ord1},
                  $$\Psi(p,0,tv)=h_1 (0,tv)+h_3 (0,tv)+h_4 (0,tv)-(n-1).$$
                  From $h_0 (0,tv)=1+\mathcal{O}(t^2 )$ we argue that
                 \begin{enumerate}
                           \item  $h_3 (0,tv)=\mathcal{O}(t^2 )$, 
                           \item  $h_1 (0,tv)=\left( -\triangle_{ _{\mathbb{S}^{n-1}}}v\right) t+\mathcal{O}(t^2 ) $,
                           \item  $h_4 (0,tv)=(n-1)-(n-1)tv+\mathcal{O}(t^2 )$,
                 \end{enumerate} 
                 and hence 
                 \begin{equation}
                        \Psi(p,0,tv)=\left( -\triangle_{ _{\mathbb{S}^{n-1}}}v-(n-1)v \right)t+\mathcal{O}(t^2 ). 
                 \end{equation}
                 It follows that
                  \begin{equation}
                             \left[\frac{\partial}{\partial x}\Psi(p,r,x))\right] _{r=0, x=0}(v)=              
                             -\triangle_{ _{\mathbb{S}^{n-1}}}v-(n-1)v.
                  \end{equation}
                  Concerning the bijectivity, the proof is an immediate application  of the following lemma  
                  with $T=L$, and this completes the proof. 
        \end{Dem}
        \begin{Lemme}\cite[Cor. 32 Appendix, page 464]{Bes}\\
              Let E,F be two complex or real vector bundles on $\mathcal{M}$.\\ 
              Let $T:C^{\infty }(E)\rightarrow C^{\infty}(F)$ be a linear differential operator of order $k$.\\
              We consider an extension of $T:L^2(E)\rightarrow L^2 (F)$.\\  
              If T is elliptic or underdetermined elliptic then 
              \begin{itemize}
                            \item $Ker(T^*)\subset C^{\infty}(F)$\\
                            \item $dim_{\mathbb{C}}(Ker(T^*))< \infty$\\
                            \item $W^{k,p}(F)=T(W^{k+l,p}(E))\bigoplus Ker(T^*)\;(1<p<+\infty) $\\
                            \item $C^{l,\alpha}(F)=T(C^{k+l,\alpha}(E))\bigoplus Ker(T^*)$\\
                            \item $C^{\infty}(F)=T(C^{\infty}(E))\bigoplus Ker(T^*)$ 
              \end{itemize}  
        \end{Lemme}  
\subsection{Differentiability of center of mass}
\begin{Lemme}
           There exists a smooth map \Fonct{c}{\mathbb{R}\times\Gamma(\mathcal{F}^{2,\alpha})}{\mathcal{M}}{(r,p,x)}{c(r,p,x)}  defined implicitly by the equation
           \begin{equation}\label{cmimpl}
                      \int_{\mathcal{N}_{p,r,x}}exp_{c}^{-1}z dVol(z)=0
           \end{equation} 
           in $T_c \mathcal{M}$ where 
$$\mathcal{N}_{p,r,x}=\{ exp_p (r(1+x(\theta))\theta)|\theta\in\mathbb{S}^{n-1} \}.$$
\end{Lemme}
\textbf{Remark:} $c(r,p,x)$ is the center of mass of $\mathcal{N}_{r,p,x}$ and $c(0,p,x)=p$.\\
 \\
\begin{Dem} We rewrite formula (\ref{cmimpl}) in the following more suitable form for our purposes:
          \begin{equation}\label{cmimpl1}
                     \int_{\mathbb{S}^{n-1}T_p \mathcal{M}}exp_{c}^{-1}(exp_p (r(1+y(p)(\theta))\theta)) f^*(dVol_{\mathcal{N}_{r,p,y(p)}})(\theta )=0 
           \end{equation} 
           in $T_c \mathcal{M}$ where 
\begin{equation}\label{perdens}
f^*(dVol_{\mathcal{N}_{r,p,y(p)}})=\sigma dVol_{(\mathbb{S}^{n-1} ,can)}
\end{equation}
\Fonct{f}{\mathbb{S}^{n-1}}{\mathcal{N}_{p,r,x}}{\theta}{exp_p (r(1+x(\theta))\theta).}
Let $o\in\mathcal{M}$ be such that $p$ and $c$ are in a normal neighborhood of $o$, we choose a local field of orthonormal frame, this choice gives an isometry $\tilde{\theta}\mapsto\theta(p,\tilde{\theta})$ of $\mathbb{S}^{n-1}T_o \mathcal{M}$ on $\mathbb{S}^{n-1}T_p \mathcal{M}$. Then (\ref{cmimpl1}) 
becomes
 \begin{equation}
                     \int_{\mathbb{S}^{n-1}T_o \mathcal{M}}exp_{c}^{-1}(exp_p (r(1+y(p)(\theta))\theta)) f^*(dVol_{\mathcal{N}_{r,y}})JdVol_{(\mathbb{S}^{n-1} ,can)}(\tilde{\theta})=0, 
           \end{equation} 
because $|det(d\theta(\tilde{\theta},p))|=J=1$.\\
We set
           \begin{equation}
                      \tilde{f}(c,r,p,x):=\int_{\mathbb{S}^{n-1}T_p \mathcal{M}}exp_{c}^{-1}(exp_p (r(1+y(p)(\theta))\theta))\sigma dVol_{\mathbb{S}^{n-1}}(\theta )
           \end{equation} 
           \Fonct{F}{\mathcal{M}\times\mathbb{R}\times\Gamma(\mathcal{F}^{2,\alpha})}{T\mathcal{M}}{(c,r,p,x)}{\tilde{f}(c,r,p,x)}
We verify that $F$ is $C^{\infty}$ with respect to all variables,
\begin{itemize}
      \item by exchanging the operation of derivation and integration, and the differentiability of the exponential map on $T\mathcal{M}$ one gets the differentiability with respect to $c$ of $F$, 
      \item $x\mapsto dx$ is smooth from $C^{2,\alpha}$ to $C^{1,\alpha}$ by continuity and linearity,
      \item $\theta\mapsto r(1+x(\theta ))\theta$ is $C^{\infty}(\mathbb{S}^{n-1},T_p \mathcal{M})$,
      \item $(r,p,x)\mapsto f^*(dVol_{\mathcal{N}_{p,r,x}})$ is $C^{\infty}(\mathbb{R}\times\mathcal{F}^{2,\alpha},C^{1,\alpha}(\mathbb{S}^{n-1},\mathbb{R}))$ because it is the norm of a multivector of $\wedge_{n-1}T_p \mathcal{M}$ whose components are the determinants in smooth functions of $r, p, x, dx$,
      \item $exp_p (r(1+x(\theta )))$ is $C^{\infty}$,
      \item we take composition with $exp_c ^{-1}$,
      \item we multiply by $\sigma$,
      \item at last we  integrate on $\mathbb{S}^{n-1}$, that is a linear continuous operation. 
\end{itemize}
Now, we can see that $F$ is divisible by $r^{n-1}$ without changing the smoothness of the resulting function $G$ and finally, that we can apply the implicit function theorem to $G=\frac{F}{r^{n-1}}$.
From the preceding arguments we see that $G$ is $C^{\infty}$ with respect to all variables.
$$\frac{\partial G}{\partial c} (c(0,p,0),0,p,0)=\frac{\partial}{\partial c}\left[\int_{\mathbb{S}^{n-1}T_p \mathcal{M}}exp_{c}^{-1}(p) dVol_{(\mathbb{S}^{n-1} ,can)}(\theta )\right]_{c=p}=\alpha_{n-1}\vartriangle $$
where $\alpha_{n-1} :=Vol_{(\mathbb{S}^{n-1},can)}(\mathbb{S}^{n-1})$ is the $(n-1)$-dimensional volume of the unit sphere with respect to the canonical metric induced by the Euclidean one of $\mathbb{R}^n$ and \Fonct{\vartriangle}{T_p \mathcal{M}}{T_p \mathcal{M}\times T_p \mathcal{M}}{y}{(y,-y).}
Let us chose a trivialisation of $T_p \mathcal{M}$ and compose $G$ with the projection $\pi_1$ on vertical fibers.
So $$\frac{\partial \pi_1 \circ G}{\partial c} (c(0,p,0),0,p,0)=\pi_1 (\alpha_{n-1}\vartriangle)=-\alpha_{n-1}Id.$$
Hence the implicit function theorem applies, there exists a unique $c(r,y)$ of class $C^{\infty}$ such that  $\pi_1 (G(c(r,y),r,y))=0$, i.e. $G(c(r,y),r,y)$ is the zero section.
\end{Dem}
\begin{Lemme}\label{cms}
            There exists a smooth map $A$ such that 
            \begin{equation}\label{lc1}
                      exp_p ^{-1}(c(r,y))=rA(r,y)
            \end{equation}
           and 
           \begin{equation}\label{lc2}
                     A(0,y)=\frac{\int_{\mathbb{S}^{n-1}}(1+x(\theta))^{n-1}\theta \sqrt{||dx||^2+(1+x)^2}}{\int_{\mathbb{S}^{n-1}}(1+x(\theta ))^{n-2}\sqrt{||dx||^2+(1+x)^2}}.
           \end{equation}
\end{Lemme}
\begin{Dem}
           In order to show (\ref{lc1}) it suffices to remark that $$c(0,y)=p.$$\\
           We now show (\ref{lc2}). We choose a trivialisation of the tangent bundle.
            \begin{Lemme}\label{exppc}
             \begin{equation}\label{ep2}
                                   exp_{c}^{-1}(exp_p (v))=exp_c ^{-1}(p)+v+o(||p-c||+||v||).
                        \end{equation}    
           This reflects the fact that Riemannian manifolds are Euclidean at small scale.
            \end{Lemme}
            \begin{Dem}
                        In an arbitrary system of coordinates of a neighborhood of $p$, (\ref{ep2}) becomes 
                        \begin{equation}
                                   exp_{c}^{-1}(exp_p (v))=p-c+v+o(||p-c||+||v||). 
                        \end{equation}
In order to show this, we set $q:=p-c$ and $H(p,q,v)=exp_{p-q}^{-1}(exp_p (v))$.\\
 Expanding the $C^{\infty}$ function    
                        $H$ at first order in a neighborhood of $(0,0)$, for all $p$, in function of $q$  and $v$. \\
                        $H(p,0,v)=v$ implies that $\frac{\partial}{\partial v}H(p,0,0)=id$.\\
                        $H(p,q,0)= exp_{p-q}^{-1} (p)=q+o(|q|)$ because if $w:= exp_{p-q}^{-1} (p)$ then 
                        $exp_{p-q}(w)=p-q+w+o(|w|)=p$ due to the fact that $d_0 exp_z =id$ for all  $z$ and hence $q=w+o(|w|)$ which implies that $w=q+o(q)$.                             
             \end{Dem}
          We choose a trivialisation of the tangent bundle and we identify $y=(p,x(p))$, then
           from the equation $G(c(r,y),r,y)=0$ we deduce $$A(r,y)=\frac{\int_{\mathbb{S}^{n-1}}(1+x(\theta))^{n-1}\theta \sqrt{||dx||^2+(1+x)^2+o(r)}}{\int_{\mathbb{S}^{n-1}}(1+x(\theta ))^{n-2}\sqrt{||dx||^2+(1+x)^2+o(r)}}+o(1).$$
           If, we put in the preceding equation $r=0$ we obtain (\ref{lc2}).\\
           This follows from the fact that 
           \begin{equation}
                       \sigma=(1+x(\theta ))^{n-2}r^{n-1}\sqrt{||dx||^2+(1+x)^2+o(r)}
           \end{equation}
           and 
           \begin{equation}
                        exp_{c}^{-1}(exp_p (r(1+y(p)(\theta))\theta))\sim exp_c ^{-1}(p)+r(1+x)\theta+o(r)
           \end{equation} 
           from which it follows
            \begin{equation}
                      \frac{1}{r^n} \int_{\mathbb{S}^{n-1}T_p \mathcal{M}} (exp_c ^{-1}(p)+r(1+x)\theta+o(r) )\sigma dVol_{(\mathbb{S}^{n-1} ,can)}=0
            \end{equation}
           and finally
           \begin{equation}
                      \int_{\mathbb{S}^{n-1}T_p \mathcal{M}}-A(r,y)\frac{\sigma}{r^{n-1}}+ \int_{\mathbb{S}^{n-1}T_p \mathcal{M}}(1+x)\frac{\sigma}{r^{n-1}}=o(1)
           \end{equation}  
           which proves the lemma.
\end{Dem}
\subsection{Existence and uniqueness of pseudo-balls}
 \begin{Lemme}\label{TFI1}
            Let  
            \Fonct{\Phi}{\mathbb{R}\times\Gamma(\mathcal{F}^{2,\alpha})}{T\mathcal{M}\times\left(  C^{0,\alpha}(\mathbb{S}^{n-1})\cap (ker(L))^{\perp} \right) }{(r,p,x)}{\left( A(r,y),Q\circ\Psi (r,p,x)\right) } 
Then, for $r$ sufficiently small,\\
there exists a  unique $x(p,r)$ in $C^{2,\alpha}(\mathbb{S}^{n-1})$ of  small $C^{2,\alpha}$ norm, solution of the implicit equation $$\Phi (r,p,x(p,r))=\Phi (r,y_r)=(y_0 ,0).$$
Here $y_0$ is the zero section of $T\mathcal{M}$.\\
Furthermore, $x$ depends smoothly on $p$ and $r$.
 \end{Lemme}
\begin{Dem}
           Like in the preceding lemma, we remark  that $\Phi$ is $C^{\infty}$ in all its arguments.
           Trivialising $T\mathcal{M}$ in  the usual manner we have:
           $$\frac{\partial \Phi}{\partial x}(p,0,0)=(\frac{d}{dt} A(0,tv)_{|t=0}, L)$$   
           \Fonct{\frac{\partial \Phi}{\partial x}(p,0,0)}{C^{2,\alpha}(\mathbb{S}^{n-1} )}{T_p\mathcal{M}\times C^{0,\alpha}(\mathbb{S}^{n-1})\cap (ker(L))^{\perp}}{v}{(\frac{n}{\alpha_{n-1}}\int_{\mathbb{S}^{n-1}}v(\theta)\theta dVol_{(\mathbb{S}^{n-1} ,can)}, L(v)).}
Remember that here $L=-\triangle_{\mathbb{S}^{n-1}} -(n-1)$.\\
In order to obtain this result we combine the lemma with the following calculations:
\begin{equation}\label{ocm}
       A(0,tv)=\frac{nt\int_{\mathbb{S}^{n-1}}v(\theta)\theta +O(t^2)}{\alpha_{n-1}+\int_{\mathbb{S}^{n-1}} (n-1)tvdVol_{can}+O(t^2)} 
\end{equation}
that follows by putting $x=tv$ in (\ref{lc2}).\\
Hence 
\begin{equation}\label{ocm1}
           \frac{d}{dt} A(0,tv)_{|t=0} =\frac{n}{\alpha_{n-1}}\int_{\mathbb{S}^{n-1}}v(\theta)\theta dVol_{(\mathbb{S}^{n-1} ,can)}. 
\end{equation}
          If we identify $T_p\mathcal{M}$ and $ker(L)$ via the isomorphism of vector spaces that maps vectors of $T_p\mathcal{M}$ to the restriction of linear functions to the sphere $\mathbb{S}^{n-1} T_p\mathcal{M}$, it can be easily seen that $\frac{\partial \Phi}{\partial x}(p,0,0)$ is an isomorphism.
From the previous discussion we conclude the existence of a smooth function $x(p,r)$ in the two variables such that
$c(r,p,x(p,r))=p$ and $H(p,r(1+x(p,r)))-\frac{n-1}{r}\in Ker(L)$. Furthermore the implicit function theorem also asserts that $x(p,r)$ is the unique small solution of these equations. 
\end{Dem}
\begin{Lemme}\label{rov}
          There exists
          $\rho:\mathcal{M}\times]r_0 ,r_0[ \rightarrow (-\rho_0 ,\rho_0 )$ of class $C^{\infty}$ defined implicitly by $$Vol_n (x(p,r))Vol_n (\mathcal{N}^+ _{p,r,x(p,r)})=\omega_n \rho(p,r)^n ,$$
and there exists $r:\mathcal{M}\times ]\rho_0 ,\rho_0 [$ defined by 
$$Vol_n (x(p,r(p,\rho ))=\omega_n \rho^n .$$
Here $\mathcal{N}^+ _{p,r,x(p,r)}=\left\lbrace exp_p (t\theta)| 0\leq t\leq r(1+x(p,r))\right\rbrace$.
\end{Lemme}
\begin{Dem}
            Because $Vol_n (x(p,r))=\omega_n r^n h(p,r)$ where $h$ is a $C^{\infty}$ function with $h(p,0)=1$, then we can write $\rho(p,r):=r h(p,r)^\frac{1}{n}$.
            Furthermore it is easy to see that $\frac{\partial\rho}{\partial r}(p,0)=1$, hence we can solve 
            $$r h(p,r)^\frac{1}{n}=\rho$$ in $r$ and obtain a function $r(p,\rho(p,r))=r$.
\end{Dem}
\begin{Res}\label{res1}
There exist $\rho_0$ and a smooth map $\beta:\mathcal{M}\times\mathbb{R}\rightarrow\mathcal{F}^{2,\alpha}$ such that for all $p\in\mathcal{M}$, and for all $\rho_0 >\rho >0$, the hypersurface $exp_p (\beta(p,\rho)(\theta )\theta) $ is the unique pseudo-ball which has its center of mass at $p$ and enclosing a volume $\omega_n |\rho|^n$.
\end{Res}
\begin{Dem}
          $\beta(p,\rho):=r(p,\rho )(1+x(p,r(p,\rho )))$ where $x(p,r)$ is produced by the preceding lemma.
\end{Dem}

      \newpage
      \section{Approximate solutions of the isoperimetric problem for small volumes are nearly round spheres}
\subsection{Introduction}
In this section it is assumed that   
\begin{enumerate}
           \item $\mathcal{M}$ has bounded geometry ($|\mathcal{K}|\leq\Lambda$ and $inj_{\mathcal{M}}\geq\varepsilon >0$) where $inj_{\mathcal{M}}$ is the injectivity radius of $\mathcal{M}$ ,
           \item the domains $D_j \in\tau_{\mathcal{M}}$ are approximate solutions i.e. $\frac{Vol_{n-1}(\partial D_j)}{I(Vol_n (D_j ))}\rightarrow 1$ for $j\rightarrow +\infty$.
\end{enumerate} 

We prove in this section the following theorem.

\begin{Thm}\label{tprec1}
Let $(\mathcal{M},g)$ be a Riemannian manifold with bounded geometry, $D_j $ a sequence of approximate solutions of the isoperimetric problem such that $Vol_g (D_j )\rightarrow 0$. Then there exist $p_j\in \mathcal{M}$, and radii $R_j$ such that  
\begin{equation}
            \lim_{j\rightarrow +\infty}\frac{Vol(D_j \Delta B(p_j , R_j ))}{Vol(D_j )}\rightarrow 0.
\end{equation}
\end{Thm}     
The proof of theorem \ref{tprec1} occupies the rest of the section.
\subsection{Taylor's theorem revisited}
Jean Taylor has shown that polyhedral chains in $\mathbb{R}^n$ which are approximate solutions of the isoperimetric problem are close to balls in the mass norm, as stated in the following theorem.
\begin{Def}
We denote by $ c_n := \frac{Vol_{n-1}(\mathbb{S}^{n-1})}{[Vol_n (\mathbb{B}^{n})]^{\frac{n-1}{n}}}$ the constant in the Euclidean isoperimetric profile.
\end{Def}
\begin{Thm}\label{JT1}
Let $W$ be the $n$-ball of $\mathbb{R}^n$ centered to the origin, of volume $1$.\\ 
           Let $\left\{ S_j \right\}\subset \mathcal{P}_n (\mathbb{R}^n )$ be a sequence of polyhedral chains (i.e. of density $1$) 
           contained in a big ball of $\mathbb{R}^n$, of barycenter at the origin, $\mathbf{M}(S_j)=1$ and  satisfying
           $$\lim_{j\rightarrow +\infty}\mathbf{M}(\partial S_j )=\mathbf{M}(\partial W ).$$
           Then  
           $$\mathbf{M}(S_j -W)\rightarrow 0.$$
\end{Thm}
\begin{Dem} Apply Taylor's theorem as stated in pages 420-421 of \cite{TJ} to the constant function $F$ equal to $1$.
\end{Dem}

It turns out that the same theorem it's true if the minimizing sequence is composed of more general currents than polyhedral chains, for example of integral currents (by the strong approximation theorem of Federer 4.2.20) and also if the minimizing sequence is not of bounded diameter. This follows from arguments which are somehow hidden in \cite{Alm}. A good reference for  the following theorem is  \cite{LeoRigot}.
\begin{Thm}\label{JT4}
           Let $\left\{ T_j \right\}\subset \mathbb{I}_n (\mathbb{R}^n )$ be a sequence of integral currents,      
           satisfying
     $$\lim_{j\rightarrow +\infty}\frac{\mathbf{M}(\partial T_j )}{\mathbf{M}(T_j )^{\frac{n-1}{n}}}=c_n .$$
           Then 
           there exist balls $W_j$ such that $$\frac{\mathbf{M}(T_j -W_j)}{\mathbf{M}(W_j )}\rightarrow 0.$$
\end{Thm}
\subsection{Lebesgue numbers}
Let $(\mathcal{M},g)$ be a Riemannian manifold with bounded geometry. We can construct a good covering of  $\mathcal{M}$ by balls having the same radius.
\begin{Lemme}\label{lebesgue}
Let $(\mathcal{M},g)$ be a Riemannian manifold with bounded geometry. There exist an integer $N$, some constants $C$, $\epsilon>0$ and a covering $\mathcal{U}$ of $\mathcal{M}$ by balls having the same radius $3\epsilon$ and having also the following properties.
\begin{enumerate}
  \item $\epsilon$ is a Lebesgue number for $\mathcal{U}$, i.e. every ball of radius $\epsilon$ is entirely contained in at least one element of $\mathcal{U}$ and meets at most $N$ elements of $\mathcal{U}$.
  \item For every ball $B$ of this covering, there exist a $C$ bi-Lipschitz diffeomorphism on an Euclidean ball of the same radius.
\end{enumerate}
\end{Lemme}

\begin{Dem} Let $\epsilon=\frac{inj_{\mathcal{M}}}{2}$. Let $\mathcal{B}=\{ B(p,\epsilon) \}$ be a maximal family of balls of $\mathcal{M}$ of radius $\epsilon$ that have the property that any pair of distinct members of $\mathcal{B} $ have empty intersection. Then the family $2\mathcal{B}:=\{ B(p,2\epsilon) \}$ is a covering of $\mathcal{M}$. Furthermore, for all $y\in\mathcal{M}$, there exist $B(p,\varepsilon )\subset\mathcal{B}$  such that $y\in B(y,2\epsilon)$ and thus $B(y,\varepsilon )\subseteq B(p,3\epsilon)$. Hence $\epsilon$ is a Lebesgue number for the covering $3\mathcal{B}$. Let $B(p,3\epsilon)$ and $B(p',3\epsilon)$ be two balls of $3\mathcal{B}$ having non empty intersection. Then $d(p,p')<6\epsilon$, hence $B(p',\epsilon)\subseteq B(p,7\epsilon)$. The ratios $Vol(B(p,7\epsilon))/Vol(B(p,\epsilon))$ are uniformly bounded because the Ricci curvature of $\mathcal{M}$ is bounded from below, and hence the Bishop-Gromov inequality applies. The number of disjoints balls of radius $\epsilon$, contained in $B(p,7\epsilon)$, is bounded and does not depend on $p$. Thus the number of balls of $3\mathcal{B}$ that intersect one of these balls is uniformly bounded by an integer $N$. We conclude the proof by taking $\mathcal{U}:=3\mathcal{B}$.
In fact by Rauch's comparison theorem, for every ball $B(p,\epsilon)$, the exponential map is 
$C$ bi-Lipschitz with a constant $C$ that depends only on $\epsilon$ and on upper bounds for the sectional curvature $\mathcal{K}$.
\end{Dem}

\subsection{Cutting domains in small diameter subdomains} 
This section is inspired by the article of B\'erard and Meyer \cite{BM} lemme II.15 and the theorem of appendix C page 531.
\begin{Prop}
Let $I$ be the isoperimetric profile of $\mathcal{M}$. Then
      $$ \limsup_{a\rightarrow 0}\frac{I(a)}{a^{\frac{n-1}{n}}}\leq c_n .$$
\end{Prop}
\begin{Dem} Fix a point $p\in\mathcal{M}$.
      \begin{eqnarray*}
            \limsup_{a\rightarrow 0}\frac{I(a)}{a^{\frac{n-1}{n}}}\leq  \limsup_{a\rightarrow 0}
            \frac{Vol(\partial B(p,r(a)))}{Vol(B(p,r(a)))^{\frac{n-1}{n}}}
      \end{eqnarray*}
      with $r(a)$ such that $Vol(B(p,r(a)))=a$.\\
      Changing variables in the limits, we find
      \begin{eqnarray*}
            \limsup_{a\rightarrow 0}\frac{Vol(\partial B(p,r(a)))}{Vol(B(p,r(a)))^{\frac{n-1}{n}}} & = & \limsup_{r\rightarrow 0}\frac{Vol(\partial B(p,r))}{Vol(B(p,r))^{\frac{n-1}{n}}}\\
         \limsup_{r\rightarrow 0}\frac{r^{n-1}Vol(\mathbb{S}^{n-1})+\cdots }{[r^{n}Vol(\mathbb{B}^{n})+\cdots ]^{\frac{n-1}{n}}}& = & c_n .      
      \end{eqnarray*}
\end{Dem}
  \begin{Def}
      Let $r>0$. We define the \emph{unit grid} of $\mathbb{R}^{n}$ and we denote by $G_1 $ the set of points which have at least one integer coordinate ($\in\mathbb{Z}$). We call \emph{grid of mesh} $r$ in  $\mathbb{R}^n$ a set $G$ of the form $v+rG_1$ where $v\in \mathbb{R}^{n}$. 
      We denote by $\mathcal{G}_r:=([0,r]^n , \mathcal{L}^n ) $ the set of all grids of mesh $r$, endowed with its natural Lebesgue mesure.       
\end{Def}
\begin{Prop}\label{l5}
        Let $D$ be an open set of $\mathbb{R}^n$.
       $$ \frac{1}{r^n}\int_{\mathcal{G}_r} Vol_{n-1}(D\cap G) \mathcal{L}^n (dG)=\frac{n}{r}Vol_n (D). $$
\end{Prop}
\begin{Dem}
      We observe that every grid $G$ decomposes as a union of $n$ sets $G^{(i)}$ of the type $v+tG_1^{(i)}$ where $G_1^{(i)}$ is the set of points with integer $i-$th coordinate.\\
Moreover $G^{(i)}\cap G^{(j)}$ has $(n-1)$-dimensional Hausdorff measure equal to zero.
      \begin{eqnarray*}
            \frac{1}{r^n}\int_{\mathcal{G}_r} Vol_{n-1}(D\cap G) \mathcal{L}^n (dG) & = & \frac{1}{r^n}\sum_{i=1} ^n  \int_{[0,r]^n }Vol_{n-1}(D\cap G^{(i)}) \mathcal{L}^n (dG)\\
   &=&         \frac{1}{r^n}\sum_{i=1} ^n  \int_{0}^r r^{n-1} Vol_{n-1}(D\cap G^{(i)}) \mathcal{L}^n (dG)\\
    & = & \frac{n}{r}Vol_n (D).
      \end{eqnarray*} 
\end{Dem}   
\begin{Cor}
Let $r>0$. Let $D$ be an open set of $\mathbb{R}^n$.  There exists a grid $G$ of mesh $r$ such that  
\begin{equation}
Vol_{n-1}(D\cap G)\leq \frac{n}{r}Vol_n (D).
\end{equation}  
\end{Cor}   
\begin{Prop}\label{l1}
      We denote $D_{G,k}$ the connected components of $D\setminus G$. Then
      $$ \frac{\sum_k Vol(\partial D_{G,k})-Vol(\partial D)}{Vol(D)^{\frac{n-1}{n}}}\rightarrow 0$$
      for $\displaystyle
          \frac{Vol(D)^{\frac{1}{n}}}{r}\rightarrow 0.$ 
\end{Prop}
\begin{Dem}
      For every grid $G$,
      $$\sum_k Vol(\partial D_{G,k})-Vol(\partial D)=2Vol_{n-1}(D\cap G).$$ By corollary \ref{l5}, there exists a grid $G$ such that $Vol_{n-1}(D\cap G)\leq \frac{n}{r}Vol_n (D)$. We deduce that 
      \begin{eqnarray*}
             0\leq\frac{\sum_k Vol(\partial D_{G,k})-Vol(\partial D)}{Vol(D)^{\frac{n-1}{n}}}\leq\frac{\frac{2n}{r}Vol_n (D)}
             {Vol(D)^{\frac{n-1}{n}}}=\frac{2nVol_n (D) ^{\frac{1}{n}}}{r}.
      \end{eqnarray*}
      Thus if $r$ is very large with respect to $Vol(D)^{\frac{1}{n}}$ then  
      $$  \frac{\sum_k Vol(\partial D_{G,k})-Vol(\partial D)}{Vol(D)^{\frac{n-1}{n}}}$$
      is close to $0$.
\end{Dem}

\begin{Prop}\label{l2}
      Let $\mathcal{M}$ be a Riemannian manifold with bounded geometry. Let $D_j$ be a sequence of domains of  $\mathcal{M}$ so that 
      \begin{enumerate}
              \item $Vol_n (D_j )\rightarrow 0$.
    \item $\limsup_{j\rightarrow +\infty} \frac{Vol_{n-1}(\partial D_j )}{Vol(D_j )^{\frac{n-1}{n}}}\leq c_n$.
      \end{enumerate}
      For any sequence $( r_j  )$ of positive real numbers that tends to zero ($r_j \rightarrow 0$ ) and  $\frac{Vol(D_j )^{\frac{1}{n}}}{r_j }\rightarrow 0$,
      there exists a cutting $D_j =\bigcup_k D_{j,k}$ of $D_j$ in domains $D_{j,k}$ with
       $Diam(D_{j,k})\leq const_{\mathcal{M}}\cdot r_j$ such that 
      $$\limsup_{j\rightarrow +\infty} \frac{\sum_k Vol_{n-1}(\partial D_{j,k})}{(\sum_k Vol(D_{j,k} ))^{\frac{n-1}{n}}}\leq c_n .$$       
\end{Prop}
\begin{Dem}
We apply lemma \ref{lebesgue} and we take a covering $\{ \mathcal{U}\}$ of $\mathcal{M}$ by balls of radius $3\epsilon$, of multiplicity $N$ and Lebesgue number $\epsilon>0$. For every ball $B(p,3\epsilon)$ of this family, we fix a diffeomorphism $\phi_p :B(p,3\epsilon)\to B_{\mathbb{R}^n}(0,3\epsilon)$ of  Lipschitz constant $C$.\\
For every $j$ we fix also a radius $r_j >>Vol_{n}(D_j)^{\frac{1}{n}}$ and we map the grids of mesh $r_j$ of $\mathbb{R}^n$ in $B(p,3\epsilon)$ via $\phi_p$, i.e. for $G\in\mathcal{G}_{r_{j}}$, we have
$$G_p =\phi_{p}^{-1}(G).$$
Let us denote by $D_{j,k}$ the connected components of $D_j \setminus (\cup_p G_p )$.
We are looking for an estimate of the supplementary volume introduced by the cutting in this $D_{j,k}$,
$$\sum_k Vol_{n-1}(\partial D_{j,k})-Vol_{n-1}(\partial D_j )=2Vol_{n-1}(D_j \cap (\cup_l G_l )).$$ 
First estimate the average $m=   \frac{1}{r_j ^n }\int_{\mathcal{G}_{r_j}}Vol_{n-1}(D_j\cap (\cup_l G_l ))\mathcal{L}^n(dG)$ of this volume over all possible choices of the grids $G\in\mathcal{G}_{r_{j}}$.
\begin{eqnarray*}
m   & \leq & 
      \frac{1}{r_j ^n }\sum_p \int_{\mathcal{G}_{r_j}}Vol_{n-1}(D_j\cap G_p )\mathcal{L}^n (dG) \\
      & \leq &  \frac{1}{r_j ^n }\sum_p \int_{\mathcal{G}_{r_j}}{Vol_{n-1}}_{(\mathbb{R}^n , {\phi_p ^{-1}}^{*} (g))}(\phi_p (D_j )\cap G)\mathcal{L}^n(dG)\\
      & \leq & 
      \frac{C}{r_j ^n }\sum_p \int_{\mathcal{G}_{r_j}}{Vol_{n-1}}_{(\mathbb{R}^n , can)}(\phi_p (D_j\cap\mathcal{U}_p  )\cap G)\mathcal{L}^n(dG)\\
      & \leq & C\frac{n}{r_j }\sum_p Vol_n (\phi_p (D_j \cap B(p,3\epsilon)))\\
      & \leq & C^2\frac{n}{r_j }\sum_p Vol_n (D_j \cap B(p,3\epsilon))\\
      & \leq & C^2 \frac{n}{r_j}NVol_{n}(D_j ).
\end{eqnarray*}
This is true because every point of $\mathcal{M}$ is contained in at most $N$ balls $B(p,3\epsilon)$. Then there exists $G $ in $\mathcal{G}_{r_j}$ such that 
$$Vol_{n-1}(D_j \cap (\cup_p G_p ))\leq C^2 \frac{n}{r_j}NVol_{n}(D_j ),$$ and so
$$0\leq\frac{\sum_k Vol_{n-1}(\partial D_{j,k})-Vol_{n-1}(\partial D_j )}{Vol_{n}(D_j )^{\frac{n-1}{n}}}\leq  2 C^2 \frac{n}{r_j}NVol_{n}(D_j )^{\frac{1}{n}}.$$
From the last inequality we obtain    
$$\limsup_{j\rightarrow +\infty} \frac{\sum_k Vol_{n-1}^{\mathcal{M}}(\partial D_{j,k})}{(\sum_k Vol_n ^{\mathcal{M}}(D_{j,k} ))^{\frac{n-1}{n}}}=
\limsup_{j\rightarrow 0} \frac{Vol_{n-1}^{\mathcal{M}}(\partial D_j )}{Vol_n ^{\mathcal{M}}(D_j )^{\frac{n-1}{n}}}\leq c_n .$$
Now, fix $x\in D_j$. By construction, $\epsilon$ is a Lebesgue number of the covering $\{ \mathcal{U}\}$, and there exists a ball $B(p,3\epsilon)$ that contains
$B_{\mathcal{M}}(x,\epsilon)$. Let $D_{j,k}$ denote the connected components of $D\setminus (\cup_p G_p )$ that contains $x$, and $D'_{j,k}$ the connected components of $\phi_p (B(p,\epsilon))\setminus G$ that contains $\phi_p (x)$. We observe that $D'_{j,k}$ is a cube of edge $r_j$, if $j$ is large enough so that $r_j \leq \epsilon/C\sqrt{n}$, then $D'_{j,k}$ is contained in $\phi_p (B(p,\epsilon))$, hence $D_{j,k}$ is contained in $\phi_{p}^{-1}D'_{j,k}$, which have diameter at most $C\,r_j$.
\end{Dem}

\subsection{Selecting a large subdomain}

We first show that an almost Euclidean isoperimetric inequality can be applied to small domains.
\begin{Lemme}\label{lem1}Let $\mathcal{M}$ be a Riemannian manifold with bounded geometry.\\
                          Then\\
                          \begin{equation}
                                   \frac{Vol(\partial D)}{Vol(D)^{\frac{n-1}{n}}}\geq c_n (1-\eta (diam(D)))
                          \end{equation} 
with $\eta\rightarrow 0$ as $diam(D)\rightarrow 0$.
\end{Lemme}
\begin{Dem}
In a ball of radius $r<inj(\mathcal{M})$, we reduce to the Euclidian isoperimetric inequality via the exponential map, that is a $C$ bi-Lipschitz diffeomorphism with $C=1+\mathcal{O}(r^2 )$. This implies for all domains of diameter $<r$, 
$$ \frac{Vol(\partial D)}{Vol(D)^{\frac{n-1}{n}}}\geq c_n C^{-2n+2}=c_n (1-\mathcal{O}(r^2 )).$$   
\end{Dem}

Second, we have a combinatorial lemma that tells us how in a cutting the largest domain contains almost all the volume.

\begin{Lemme}\label{l3}
      Let $f_{j,k}\in [0,1]$ be numbers such that for all $j $, $\sum_k f_{j,k} =1$.
      Then
      $$ \limsup_{j\rightarrow + \infty} \sum_k f_{j,k}^{\frac{n-1}{n}}\leq 1 $$
      implies that
      $$ \lim_{j\rightarrow +\infty } \max_k  f_{j,k}=1.$$
\end{Lemme}
\begin{Dem}
      We argue by contradiction.
      Suppose there exists $\varepsilon >0$ for which there exists $j_{\varepsilon}\in \mathbb{N}$ so that for all  $j\geq j_{\varepsilon}$, we have $\max_k \{f_{j,k}\}\leq 1-\varepsilon $.
      Then for all $j\geq j_{\varepsilon}$, we have $f_{j,k}\leq 1-\varepsilon $. From this inequality,
 $$\sum_k f_{j,k}^{\frac{n-1}{n}}=\sum_k f_{j,k}f_{j,k}^{\frac{-1}{n}}\geq \frac{\sum_k f_{j,k}}{(1-\varepsilon )^{\frac{1}{n}}}
        \geq \frac{1}{(1-\varepsilon)^{\frac{1}{n}}},$$ 
      hence
      $$\limsup_{j\rightarrow + \infty} \sum_k f_{j,k}^{\frac{n-1}{n}}\geq \frac{1}{(1-\varepsilon)^{\frac{1}{n}}}> 1,$$
      which is a contradiction.  
\end{Dem}
 \begin{Prop}\label{l4}
Let $\mathcal{M}$ be a Riemannian manifold with bounded geometry. Let $D_j$ be a sequence of approximate solutions in $\mathcal{M}$ with volumes that tend to zero. Let $r_j$ be a sequence of positive real numbers such that $r_j \rightarrow 0$ and $\frac{Vol(D_j )^{\frac{1}{n}}}{r_j }\rightarrow 0$.\\
There exist $p_j \in\mathcal{M}$ and $\varepsilon_j \leq const_{\mathcal{M}}r_j$ and subdomains $D'_j\subset D_j$ such that 
\begin{enumerate}
            \item $D'_j \subseteq B(p_j ,\varepsilon_j )$
            \item $\frac{Vol(\partial D'_j)}{Vol(D'_j)^{\frac{n-1}{n}}}\rightarrow 0$
            \item $\lim_{j\rightarrow +\infty } \frac{{Vol_n }^{\mathcal{M}}(D_j )}{{Vol_n }^{\mathcal{M}}(D_j )}=1.$
\end{enumerate} 
 \end{Prop}
 \begin{Dem}
       Apply proposition \ref{l2}. By definition of isoperimetric profile and lemma \ref{lem1} we have\\
       $$Vol(\partial D_{j,k})\geq I(Vol(D_{j,k} ))\geq c_n 
Vol(D_{j,k})^{\frac{n-1}{n}}(1-\eta_j )$$ where $\eta_j \rightarrow 0$. Since
  \begin{eqnarray*}
       \limsup_{j\rightarrow +\infty }\frac{\sum_k c_n 
Vol(D_{j,k})^{\frac{n-1}{n}}(1-\eta_j )}{Vol(D_j )^{\frac{n-1}{n}}} \leq
       \limsup_{j\rightarrow +\infty }\frac{\sum_k Vol(\partial 
D_{j,k})}{Vol(D_j )^{\frac{n-1}{n}}} \leq c_n ,\\ 
  \end{eqnarray*}
 $$\limsup_{j\rightarrow +\infty }\frac{\sum_k 
Vol(D_{j,k})^{\frac{n-1}{n}}}{Vol(D_j )^{\frac{n-1}{n}}}\leq \limsup_{j\rightarrow +\infty }
   \frac{1}{1-\eta_j }=1.$$
 Now, we set $f_{j,k}=\frac{Vol(D_{j,k})}{Vol(D_j)}$. We can suppose that $f_{j,1}=max_k \{f_{j,k}\}$. We apply lemma 
\ref{l3} and we deduce that  
 $$\frac{Vol(D_{j,1})}{Vol(D_j )}\rightarrow 1.$$
 But by construction $D_{j,1}\subset B_{\mathcal{M}}(p_j , 
const_{\mathcal{M}}r_j )$ with $(p_j )$ sequence of points $p_j$ in $\mathcal{M}$.\\
Finally, proposition \ref{l2} gives 
$$\limsup \frac{Vol(\partial D_j,1 )}{Vol(D_j )^{\frac{n-1}{n}}}\leq\limsup\frac{}{}\leq c_n .$$
Thus one can take $D'_j =D_{j,1}$
 \end{Dem}
\subsection{Application of Taylor's Theorem}
Let $D_j$ be a sequence of approximate solutions with $Vol(D_j )\rightarrow 0$.
According to proposition \ref{l4} there exist subdomains $D'_j \subseteq D_j$, points $p_j \in \mathcal{M}$ and radii $\varepsilon_j \rightarrow 0$ such that 
\begin{enumerate}
           \item $D'_j \subseteq B(p_j ,\varepsilon_j )$.
           \item $\frac{Vol(D'_j )}{Vol( D_j )}\rightarrow 0$.
           \item $\frac{Vol(\partial D'_j)}{Vol(D_j)^{\frac{n-1}{n}}}$.
\end{enumerate} 
We identify all tangent spaces $T_{p_j }\mathcal{M}$ with a fixed Euclidean space $\mathbb{R}^n$ and consider the domains $D''_j = exp^{-1}(D'_j )$ in $\mathbb{R}^n$. Since the pulled back metrics $\tilde{g}_j=exp_{p_j}^{*}(g_{\mathcal{M}})$ converge to the Euclidean metric, 
$$\frac{Vol(\partial D''_j)}{Vol(D''_j)^{\frac{n-1}{n}}}\rightarrow c_n. $$
According to theorem \ref{JT4}, there exist Euclidean balls $W_j =B_{eucl.}(\tilde{q}_j , R_j)$ in $\mathbb{R}^n$ such that 
$$\frac{Vol_{eucl.}(D''_j\Delta W_j)}{Vol_{eucl.}(D''_j)}\rightarrow 0.$$
Note that $\tilde {g}_j $-balls are close to Euclidean balls,
$$\frac{Vol_{eucl.}(D''_j\Delta W_j)}{Vol_{eucl.}(W_j)}\rightarrow 0.$$
Thus 
$$\frac{Vol_{eucl.}(D''_j \Delta B^{\tilde{g}_j}(\tilde{q}_j ,R_j))}{Vol_{eucl.}(D''_j)}\rightarrow 0,$$
and then, for $q_j =exp_{p_j}(\tilde{q}_j)$,
$$\frac{Vol_{eucl.}(D'_j \Delta B^{g}(\tilde{q}_j ,R_j))}{Vol_{eucl.}(D'_j)}=\frac{Vol_{eucl.}(D''_j \Delta B^{\tilde{g}_j}(\tilde{q}_j ,R_j))}{Vol_{\tilde{g}}(W_j)}\rightarrow 0.$$
Finally, since $\frac{Vol{D_j \Delta D'_j}}{Vol(D_j)}\rightarrow 0$, 
$\frac{Vol_{g}(D_j\Delta B(q_j , R_j))}{Vol_{g}(D_j)}\rightarrow 0.$ \\
This completes the proof of theorem \ref{tprec1}  
\subsection{Case of exact solutions}

\textbf{Remark:} When we consider the \emph{solutions} of the isoperimetric problem (this is the case treated in \cite{MJ}), and not \emph{approximate solutions}, the conclusion is stronger. In fact we can prove directly by the monotonicity formula that $D_j$ is of smal diameter and we can apply Taylor's theorem to the dilated of $D_j$ without passing throught any kind of cutting procedure.\\

\begin{Lemme}
      Assume $D_j$ is asolution of the isoperimetric problem. The dilated domains $D'''_j :=\frac{exp_{p_j}^{-1}(D_j )}{Vol_g (D_j)^{\frac{1}{n}}}$ are of bounded diameter and hence we can find a positive constant $R>0$ in the proof of the preceding theorem so that for all $j\in \mathbb{N}$ we have $$D'''_j \subseteq B(0,R).$$ 
\end{Lemme}
\begin{Dem}
      For the domains $D'''_j$, the mean curvature of the boundary in 
      $(\mathbb{R}^n,eucl)$ $h_j^{eucl}\leq M=const.$ for all $j$ (apply the L\'evy-Gromov isoperimetric inequality \cite{Gr1}, \cite{Gr2}) and  hence the monotonicity formula of \cite{All}[5.1 (3)] page 446 gives for a fixed $r_0$ and all $j$
      \begin{equation}
      ||\partial D'''_j ||(B(a_j , r_0))\geq e^{-M r_0}\Theta^{n-1}(||\partial D'''_j||, a_j)\omega_{n-1}r_0^{n-1}
      \end{equation}
      $a_j\in spt||\partial D'''_j ||$, $r_0$ for a fixed $r_0$ and all $j$.
      We argue 
      $$const\geq Vol_{g_{can},n-1} (\partial D'''_j )\geq \left[ \frac{Diam_{g_{can}}(D'''_j )}{2r_0}\right]\omega_{n-1}r_0 ^{n-1}$$      
      and we can conclude that $Diam_{g_{can}}(D'''_j )$ are uniformly bounded.
\end{Dem}

      \newpage
      \section{Application to compact Riemannian manifolds}

\subsection{From the pseudo-balls viewpoint}

 In this section, $\mathcal{M}$ is a \emph{compact} Riemannian manifold.
 
\begin{Lemme}\label{tprec}
Let $D_j$ be a sequence of solutions of the isoperimetric problem such that $Vol_g (D_j )\rightarrow 0$. Then eventually extracting a subsequence, there exist a point $p\in \mathcal{M}$ such that   
the domains $D_j $ are graphs in polar normal coordinates of funcions $u_j $ of class $C^{2,\alpha }$ on the unit sphere of $T_p \mathcal{M}$ of the form  
$u_j =r_j(1+v_j )$ with 
$||v_j ||_{C^{2,\alpha }(\partial B_{T_p \mathcal{M}}(0, R'))}\rightarrow 0$ and radii $r_j$.
\end{Lemme}
\begin{Dem}
Theorem \ref{tprec1} provides points $p_j$ and radii $r_j$ such that $D_j$ is close to $B(p_j , r_j )$ volumewise. Since $\mathcal{M}$ is compact, one can assume that $p_j$ converges to $p$. Let $T_j$ be $exp_p^{-1}(D_j)$ rescaled by $\frac{1}{r_j}$. Then $T_j$ is a solution of the isoperimetric problem for the rescaled pulled back metric $g_j=\frac{1}{r_j ^2}exp_p^* (g)$ which converges volumewise to a unit ball. Since the sequence $g_j$ converges smoothly to a Euclidean metric, the regularity theorem of \cite{Nar} applies hence $\partial T_j$ is the graph in normal coordinates of a smooth function $v_j$ on the unit sphere and $||v_j ||_{C^{2,\alpha}(\mathbb{S}^{n-1})}$ tends to zero.
In other words, $\partial D_j$ is the normal graph of $u_j=r_j(1+v_j )$.  
\end{Dem}
We can rewrite this lemma in the following form.
\begin{Lemme}
Let $\mathcal{M}$ be a compact Riemannian manifold of class $C^3 $.
There exists $v_1 >0$ so that for all current $D$ solution of the isoperimetric problem with  $M(D)\leq v_1 $, there exists a point $p_{D}\in\mathcal{M}$ (depending on $D$) such that $D$ is the normal graph of a function $u_D \in C^{2,\alpha }(\mathbb{S}^{n-1}) $ with $u_{D}=r_{D}(1+v_{D})$ and $||v_{D}||_{C^{2,\alpha }(\mathbb{S}^{n-1})}\rightarrow 0$ when $ M(D)\rightarrow 0$. 
\end{Lemme}
\begin{Dem}
By contradiction using the preceding lemma. \\
 If the thesis were not true then there exist $p\in \mathcal{M}$, a sequence $D_j $ of solutions of the isoperimetric problem with volumes 
  $Vol(D_j)\rightarrow 0$, and for which $\partial D_j$ is not the graph on the sphere  
$\mathbb{S}^{n-1}$ of $T_p \mathcal{M}$ of a function $u_j =r_j (1+v_j )$ where $||v_j ||_{C^{2,\alpha}}$ goes to $0$. This is in contradiction with the preceding lemma.
\end{Dem} 

\begin{Res}\label{tf} There exist $v_0$ such that if $v<v_0$ then all current solutions of the isoperimetric problem with volume $v$ are pseudo-balls.  
\end{Res}
\begin{Dem}
Let $v_1$ be as in the preceding lemma and $\rho_0$ be given by lemma \ref{rov}. We set 
$v_0 :=\min\left\lbrace \omega_n \rho_0 ^n , v_1 \right\rbrace$.
Assume $v<v_0 $. Let $T$ be a current solution of the isoperimetric problem with volume $v$. Then $\partial T$ is the normal graph centered at a point $p$ of a function $u\in C^{2,\alpha}(\mathbb{S}^{n-1})$,  $u=r(1+v)$ with $||v ||_{C^{2,\alpha}}$ small and $r$ small.
Let $q$ be the center of mass of $\partial T$. According to lemma \ref{cms}, $exp_p^{-1}(q)=rA(r,p,v)$ where $A$ is smooth and $A(0,p,0)=0$. Thus $d(p,q)\leq const(r^2+r||v||_{C^1})\leq\varepsilon r$ for $\varepsilon$ arbitrarily small. It follows that the radial projection of $ \partial T$ onto the sphere $ \partial B(q,r)$ is $C^{\infty}$ close to identity. As a consequence, $\partial T $ is the normal graph centered at q of a function $\tilde{u}=r(1+\tilde{v})$ on $T^1_q \mathcal{M}$, with $||\tilde{v}||_{C^{2,\alpha}}$ small.
Since $q$ is the center of mass of $\partial T$, $\int_{\mathbb{S}^{n-1}}\tilde{v}(\theta )\theta d\theta=0$, i.e. $\tilde{v}$ belongs to the space $C^{2,\alpha}_1(\mathbb{S}^{n-1}) $. Since $\partial T$ has constant mean curvature, it satisfies $Q(\Psi(r,q,\tilde{v}))=0$. Therefore $\partial T$ coincides with the pseudo-ball $\beta(q,\tilde{r})$ where $Vol(T)=\omega_n \tilde{r}^n$.  
\end{Dem}

\begin{Cor}
           Let $T$ be a solution of the isoperimetric problem with small enclosed volume $v$, let  $p\in\mathcal{M}$ be its center of mass. Let $St_p \leq Isom(\mathcal{M})$ be the stabilizer of $p$ for the canonical action of the group of isometries $Isom(\mathcal{M})$ of $\mathcal{M}$.\\ 
Then for all $k\in St_p$, $k(T)=T$.
\end{Cor}
\begin{Dem} Following theorem \ref{res1}, $\partial T$ is the pseudo-ball $\beta(p,r)$ where $\omega_n \rho^n =Vol_n (T)$.\\
If $k\in St_p$, then, $k(\beta(p,r))=\beta(k(p),r)=\beta(p,r)$ hence $k(T)=T$.
\end{Dem}

\subsection{Asymptotic expansion of the isoperimetric profile}

In the preceding section, we reduced the variational problem with volume constraint $v$ smaller than $v_0$ to an optimization problem on the set of pseudo-balls of enclosed volume $v$.\\
It is natural at this moment to consider the function 
\Fonct{f}{\mathcal{M}\times ]0,v_0 [}{[0,+\infty[}{(p,v)}{Vol_{n-1}(\mathcal{N}_{p,r(p,\rho)})}
where $\mathcal{N}_{p,r}=\{ exp_p \left( r(1+x(p,r)(\theta ))\right) |\:\theta\in T^1_p\mathcal{M}\}$ is the pseudo-ball of center of mass $p$, of enclosed volume $v=\omega_n \rho^n$. The following result is a reformulation of theorem \ref{tf}.

\begin{Prop} For all $v<v_0$ it is true that
           $$I_{\mathcal{M}}(v)=Inf_{p\in\mathcal{M}}\{f(p,v)\}.$$
\end{Prop}

We give an asymptotic expansion of the function $v\mapsto f(p,v)$. We use unpublished results of   Pacard and Xu. For completeness sake, the proof of the following theorem is included. Furthermore, we agree that any term denoted $\mathcal{O} ({r^k})$ is a smooth function on $\mathbb{S}^{n-1}$ that might depend on $p$ but which is bounded by a constant independent of $p$ times $r^k$ in $C^2 $ topology.

\begin{Thm}\label{PX} \emph{Asymptotic expansion of the area of pseudo-balls with respect to the radius $r$ of perturbed geodesic spheres}. 
\begin{equation}\label{daexp1}
                Vol_{n-1}(\mathcal{N}_{p,r})=r^{n-1}\alpha_{n-1}\left( 1-\frac{1}{2n}Sc(p)r^2 +\mathcal{O} ({r^4})\right). 
\end{equation}
\begin{equation}\label{daexp}
           Vol_n (\mathcal{N}_{p,r}^+ )= \frac{r^n}{n}\alpha_{n-1}\left( 1+\gamma_n (p)r^2+\mathcal{O} (r^4 )\right) 
\end{equation} 
with $\gamma_n :=-\frac{n+1}{2(n-1)(n+2)}Sc(p)$\\ and
$\mathcal{N}_{p,r}^+ :=\{ exp_p (t\theta) | 0\leq t\leq r(1+x(p,r)(\theta))\:,\:\theta\in T^1_p\mathcal{M}\} $.
\end{Thm}   
\begin{Dem}
To prove (\ref{daexp1}) and  (\ref{daexp}) we need some preliminary lemmas that allows one to expand the required quantities highlightening the geometrical meaning of the coefficients of the respective asymptotic expansions.
\begin{Lemme} Asymptotic expansion of the outward mean curvature of geodesic spheres of radius $r$.
          \begin{equation}\label{Hexp} 
                H^r_{\theta}(r,\theta)=-\frac{(n-1)}{r}+\frac{1}{3}Ric(p)(\theta )r+\mathcal{O} (r^2)
          \end{equation}
\end{Lemme}
\begin{Dem}
           Denote by $U(r,\theta)$ the shape operator of geodesic spheres of radius $r$ in the direction $\theta$ considered, as usual, as a linear operator on a fixed finite dimensional real vector space, $T_p \mathcal{M}$ for example. By standard results in Riemannian geometry (see \cite{Chav}) $U$ satisfy the following Riccati equation  
\begin{equation}\label{Riccati}
                  U'+U^2+R=0
\end{equation}
where $R$ is a suitable curvature operator and primes means derivatives taken with respect to the $r$ variable.  We are looking for an asymptotic expansion of $H=tr(-U)$.
To this aim, we start by observing that $U=J'J^{-1}$ for $J$ being the matrix whose entries are the components of the Jacobi fields vanishing at the origin, with respect to a parallel transported orthonormal basis of $T_p \mathcal{M}$ and with the initial condition $J'(0)=I$.
This means that $ \frac{J(r)}{r}\rightarrow I$ when $r\rightarrow 0$ what implies $rJ^{-1}\rightarrow I$ when $r\rightarrow 0$. For this reason $rU(r)\rightarrow I$ when $r\rightarrow 0$.
The last argument allows us to have the following asymptotic expansions
 $$U=\frac{1}{r}I+A+Br+\cdots$$
hence $$U^2 = \frac{1}{r^2}I+\frac{2}{r}A+A^2+2B+\cdots$$ 
and
$$U'=-\frac{1}{r^2}I +B+\cdots$$
where $A$ and $B$ are unknown linear operators.
Putting all these expansions in the Riccati equation (\ref{Riccati}) we obtain
the conditions on $A$ and $B$. So $A=0$ and $B=-\frac{1}{3}R$.\\ 
Finally taking traces we get equation (\ref{Hexp}), the desired result.
\end{Dem}

\begin{Lemme}
                         Let $x(p,r)$ as in lemma \ref{TFI1} then 
                         \begin{equation}
                                       x(p,r)(\theta)=r^2 x_2(p,\theta )+r^3 x_3 (p,r,\theta)
                         \end{equation}
with $x_2(p,\theta)\in C_1^{2,\alpha}(\mathbb{S}^{n-1})$, and
\begin{equation}\label{firstterm}
Lx_2=\frac{1}{3}Ric(p).
\end{equation}
\end{Lemme} 
\begin{Dem}
We already know that $x$ has an asymptotic expansion without terms of zero degree in $r$.
Let $x(r,p,\theta)=rx_1(p,\theta) +r^2 x_2 (p,\theta)+\mathcal{O}(r^3)$, we put it into 
\begin{equation}\label{implicita}
      \left( A(r,x(p,r)), Q\circ\Psi (r,p,x)\right) =(0,0),
\end{equation} 
and combining with lemma \ref{lc11} and equation (\ref{ocm1}) we must have $$Lx_1 =0,$$ and 
$$\frac{n}{\alpha_{n-1}}\int_{\mathbb{S}^{n-1}}x_1(\theta)\theta dVol_{(\mathbb{S}^{n-1} ,can)}=0,$$
the last equation means that $x_1\in C^{2,\alpha}_1$ hence $x_1\in C^{2,\alpha}_1 \cap Ker(L)=\{ 0\}$ (i.e. $x_1=0$). 
At this point, the more detailed expansion of $x$ gives $$x(r,p,\theta)=r^2 x_2 (p,\theta)+\mathcal{O}(r^3) $$ and in the same manner, we substitute analogously to what is already done in this proof we put this latter equation into  (\ref{implicita}), this yields $$Lx_2 =\frac{1}{3}Ric(p),$$ because 
we have the following expansion 
 $$rH(r(1+x))=(n-1)+Lx+r^2 h(r,p,x)$$ for a smooth function $h(r,p,x)=h(r,p,0)+\tilde{h}(r,p,x)$ with
$\tilde{h}(r,p,x)$ containing in his Taylor expansion terms that vanish at least linearly in $x$ at $x=0$.
To show (\ref{firstterm}) we observe that $\tilde{h}(r,p,r^2 x_2+\cdots)=r^2 \tilde{h}'(r,p,x_2+\cdots)=\mathcal{O}(r^2)$ hence we get
$$Q\left( Lx_2+h(0,p,0)\right) =0,$$ by equating to zero the coefficient of $r^2$ of the resulting asymptotic expansion, but $h(0,p,0)=\frac{\partial^2 \Psi}{\partial r^2}(0,p,0)=-\frac{\partial^2 H^r_{\theta}(r(1+0),\theta)}{\partial r^2}_{|r=0}=-\frac{1}{3}Ric(p)$ and this is easy to see by differentiating $rH(r(1+0))=h_4(r,0)$ twice with respect to $r$ (see calculations of proposition \ref{ord1}).
This implies (\ref{firstterm}) by observing
  $Q\left( \frac{1}{3}Ric(p)\right)=\frac{1}{3}Ric(p)$, ($Ric(p)$ is the restriction to $\mathbb{S}^{n-1} $ of a quadratic form on $\mathbb{R}^n$) and $Q(Lx_2) = Lx_2$ by definition of $L$ and $Q$.\\
We finally observe that it must be
$$\frac{n}{\alpha_{n-1}}\int_{\mathbb{S}^{n-1}}x_2(\theta)\theta dVol_{(\mathbb{S}^{n-1} ,can)}=0,$$
by equating to zero the coefficient of $r^2$ in the expansion of $A(r,r^2x_2(p)+\cdots)$ and this proves that $x_2(p,\theta)\in C_1^{2,\alpha}(\mathbb{S}^{n-1})$.\\
This finishes the proof of the lemma.
\end{Dem}

Now we need an asymptotic expansion of the $(n-1)$-dimensional volume of perturbed normal graphs on geodesic spheres.
\begin{Lemme}Let $\sigma$ be defined by equation (\ref{perdens}) then
                $$\sigma(r,p,\theta)= r^{n-1} \left[ 1+(n-1)r^2 x_2 -\frac{1}{6}Ric_p(\theta)r^2 +\mathcal{O} (r^3)\right]. $$
\end{Lemme}
\begin{Dem}
           We proceed by computing first an expansion of the pulled-back metric $g_{ij}$ on the unit sphere from $\mathcal{N}_{p,r,x}$ the general perturbed normal graph, then we compute $\sigma=\sqrt{det(g_{ij})}$.\\
Set $t=r(1+x)$
\begin{equation}
      \begin{array}{ccc}
           g_{ij} & = & <d(exp_p)_{|t}(rdx(e_i )\theta +te_i ),d(exp_p)_{|t}(rdx(e_j )\theta +te_j )>\\
                     & = & t^2 \left(\delta_{ij}+\frac{dx(e_i)dx(e_j)}{(1+x)^2}-\frac{1}{3}<R_{\theta ,e_i}\theta ,e_j>t^2 + t^4 (\cdots) \right), 
      \end{array}
\end{equation}
\begin{equation}
      \begin{array}{ccc}
           \sigma & = & t^{n-1}\sqrt{det\left(\delta_{ij}+\frac{dx(e_i)dx(e_j)}{(1+x)^2}-\frac{1}{3}<R_{\theta ,e_i}\theta ,e_j>t^2 + t^4 (\cdots) \right) }.
      \end{array}
\end{equation}
Now we put in the last equation the expansion $x=r^2 x_2+\cdots$
\begin{equation}
      \begin{array}{ccc}
           \sigma & = & r^{n-1}\left(1+(n-1)x_2 r^2+\cdots \right) \left(1-\frac{1}{6}Ric_p (\theta)+r^4(\cdots) \right) \\
                       & = & r^{n-1}\left\lbrace 1+\left[ (n-1)x_2-\frac{1}{6}Ric_p (\theta)\right] r^2\right\rbrace  +\cdots\\
      \end{array}
\end{equation}
\end{Dem}

We now come back to the proof of the theorem.\\
Integration of (\ref{firstterm}) over $\mathbb{S}^{n-1}$ yields 
$$(n-1)\int_{\mathbb{S}^{n-1}}x_2 dVol_{(\mathbb{S}^{n-1},can)}=-\frac{1}{3}\int_{\mathbb{S}^{n-1}}Ric(p) dVol_{(\mathbb{S}^{n-1},can)}.$$
To complete the proof of (\ref{daexp1}) observe that 
$$Vol_{n-1}(\mathcal{N}_{p,r})=\int_{\mathbb{S}^{n-1}}\sigma dVol_{(\mathbb{S}^{n-1},can)}$$ 
and by using the fact that 
$$\int_{\mathbb{S}^{n-1}}Ric(p) dVol_{(\mathbb{S}^{n-1},can)}=\frac{1}{n}Vol_{(n-1,can)}(\mathbb{S}^{n-1})Sc(p),$$
we get (\ref{daexp1}).\\
The proof of (\ref{daexp}) is easier and requires only the expansion of the volume density in normal polar coordinates as a function of the distance to the origin.
\begin{equation*}
           \begin{array}{lll}
                       Vol_n (\mathcal{N}_{p,r}^+ ) & = & \int_{\mathbb{S}^{n-1}}\left\lbrace \int_0 ^{r(1+x)}s^{n-1}\left[ 1-\frac{1}{6}Ric(p)s^2 \right] ds\right\rbrace dVol_{(\mathbb{S}^{n-1},can)}+\cdots\\
   & = & \int_{\mathbb{S}^{n-1}}\left\lbrace \frac{r^n}{n}+r^{n+2}\left[ x_2 -\frac{1}{6(n+2)}Ric(p) \right]\right\rbrace dVol_{(\mathbb{S}^{n-1},can)}+\cdots\\
          & = & \frac{\alpha_{n-1}r^n}{n}\left\lbrace  1- \frac{r^2}{3}\left[ \frac{1}{n-1}+\frac{1}{2(n+2)}\right] Sc(p)\right\rbrace +\cdots\\
   & = & \frac{\alpha_{n-1}r^n}{n}\left\lbrace  1- \frac{n+1}{2(n-1)(n+2)}Sc(p)r^2\right\rbrace +\cdots\\
           \end{array}
\end{equation*}
$Vol(\mathcal{N})$ and $Vol(\mathcal{N}^+)$ are even functions of $r$ variables, hence there are no terms in $r^3$  in their respective asymptotic expansions. 
By compactness of $\mathcal{M}$, the remainders are uniformly bounded in $p$.
These remarks and the last equation complete the proof. 
\end{Dem} 
   
\begin{Lemme}\label{PXvol} Asymptotic expansion of the area of pseudo-balls as a function of the enclosed volume.  
 \begin{equation}\label{Ipexp}
f(p,v)=c_n v^{\frac{n-1}{n}}\left\lbrace 1+a_p \left(\frac{v}{\omega_n} \right) ^{\frac{2}{n}} +\mathcal{O} (v^{\frac{4}{n}}) \right\rbrace  
\end{equation} 
with $a_p:=-\frac{1}{2n(n+2)}Sc(p)$.                        
\end{Lemme}

\begin{Dem}
We set $v=Vol_n (\mathcal{N}_{p,r}^+ )$ so
\begin{equation}
           v= \frac{r^n}{n}\alpha_{n-1}\left( 1+\gamma_n (p)r^2+\mathcal{O} (r^4 )\right) 
\end{equation}
We reverse the latter asymptotic expansion to obtain an asymptotic expansion of $r$ as a function of $v$. Then we substitute this expansion in (\ref{daexp1}) and we get equation (\ref{Ipexp}). 
In fact,
\begin{equation}
r^n =\frac{n}{\alpha_{n-1}}v-\gamma_n \left( \frac{n}{\alpha_{n-1}}\right)^{\frac{n+2}{n}} +\cdots
\end{equation}
\begin{equation}
r^{n-1} =\left( \frac{n}{\alpha_{n-1}}\right) ^{\frac{n-1}{n}}v ^{\frac{n-1}{n}}\left[ 1-\frac{n-1}{n}\gamma_n \left( \frac{n}{\alpha_{n-1}}\right)^{\frac{2}{n}}v^{\frac{2}{n}}+\cdots\right] 
\end{equation}
\begin{equation}
r^{n+1} =\left( \frac{n}{\alpha_{n-1}}\right) ^{\frac{n+1}{n}}v ^{\frac{n+1}{n}}\left[ 1-\frac{n+1}{n}\gamma_n \left( \frac{n}{\alpha_{n-1}}\right)^{\frac{2}{n}}v^{\frac{2}{n}}+\cdots\right] 
\end{equation}
\begin{equation}
\begin{array}{lll}
f & = & \alpha_{n-1}r^{n-1}-\frac{\alpha_{n-1}}{2n}Sc(p)r^{n+1}+\cdots\\
  & = & \alpha_{n-1}r^{n-1}-\frac{\alpha_{n-1}}{2n}Sc(p)\left( \frac{n}{\alpha_{n-1}}\right) ^{\frac{n+1}{n}}v ^{\frac{n+1}{n}}+\cdots\\
  & = & c_n v ^{\frac{n-1}{n}}\left[1-\left( \frac{n}{\alpha_{n-1}}\right)^{\frac{2}{n}}\left(\frac{n-1}{n}\gamma_n + \frac{1}{2n}Sc(p)\right) v^{\frac{2}{n}}\right] +\cdots\\
   & = & c_n v ^{\frac{n-1}{n}}\left[1-\frac{1}{2n(n+2)}Sc(p)\left( \frac{n}{\alpha_{n-1}}\right)^{\frac{2}{n}}v^{\frac{2}{n}}\right] +\cdots\\ 
\end{array}
\end{equation}
here it is no difficult to check that $\omega_n =\frac{\alpha_{n-1}}{n}$.\\
Finally it follows that 
\begin{equation*}
f(p,v)=c_n v^{\frac{n-1}{n}}\left\lbrace 1+a_p \left(\frac{v}{\omega_n} \right) ^{\frac{2}{n}} +\mathcal{O}_p (v^{\frac{4}{n}}) \right\rbrace . 
\end{equation*}  
\end{Dem}

\begin{Lemme}\label{maxsc}
           Let $v_k$ be a sequence of volumes tending to $0$. Let $T_k$ be a solution of the isoperimetric problem with enclosed volume $v_k$. Let $p_k \in\mathcal{M}$ be the center of mass of $\partial T_k$. Suppose that the sequence $p_k$ converges to a point $p\in\mathcal{M}$.\\ Then $p$ is a point of global maximum of the scalar curvature function of $\mathcal{M}$ i.e. $Sup_{\tilde{p}\in\mathcal{M}}\{ Sc(\tilde{p})\}=Sc(p)$. Furthermore,
\begin{eqnarray}\label{dl1}
\lim_{k\to\infty} I_{\mathcal{M}}(v_k )v_{k}^{\frac{-1-n}{n}}-c_n v_{k}^{\frac{-2}{n}}=\frac{c_n}{\omega_n ^{\frac{2}{n}}}a_p.
\end{eqnarray}
\end{Lemme}

\begin{Dem}
By definition of the isoperimetric profile,
$$I_{\mathcal{M}}(v_k )=f(p_k , v_k )=Inf_{p\in\mathcal{M}} \{f(p, v_k )\}.$$
           We consider the function $g(\tilde{p},v)=\left[ \frac{f(\tilde{p}, v)}{v^{\frac{n-1}{n}}}-c_n \right]\frac{1}{v^{\frac{2}{n}}}$.\\
It is easy to see that $f(p_k , v_k )=Inf_{\tilde{p}\in\mathcal{M}} \{f(\tilde{p}, v_k )\}$ if and only if $g(p_k , v_k )=Inf_{\tilde{p}\in\mathcal{M}} \{g(\tilde{p}, v_k )\}$. By lemma \ref{Ipexp}, we know that  $g(p,v)$ tends to  uniformly in $p$.has the following expansion in a neighborhood of $0$.
$$\frac{c_n}{\omega_n ^{\frac{2}{n}}}a_p$$ as $v$ tends to $0$, uniformly in $p$.\\
It follows that
\begin{itemize}
           \item $g(p_k, v_k )=Inf_{\tilde{p}\in\mathcal{M}} \{g(\tilde{p}, v_k )\}\rightarrow \frac{c_n}{\omega_n ^{\frac{2}{n}}}a_p$, 
           \item $ Inf_{\tilde{p}\in\mathcal{M}} \{g(\tilde{p}, v_k )\}\rightarrow Inf_{\tilde{p}\in\mathcal{M}} \{g(\tilde{p}, 0)\}=Inf_{\tilde{p}\in\mathcal{M}}\left\lbrace \frac{c_n}{\omega_n ^{\frac{2}{n}}}a_{\tilde{p}}\right\rbrace$. 
\end{itemize}
 hence $Inf_{\tilde{p}\in\mathcal{M}}\left\lbrace \frac{c_n}{\omega_n ^{\frac{2}{n}}}a_{\tilde{p}}\right\rbrace =\frac{c_n}{\omega_n ^{\frac{2}{n}}}a_p$ and by the presence of minus sign $-$ in the coefficient  $a_p$, we can conclude that $Sc(p)=Sup_{\tilde{p}\in\mathcal{M}}\{ Sc(\tilde{p})\}$.\\
To show (\ref{dl1}) it is sufficient to observe that 
$$I_{\mathcal{M}}(v_k )v_{k}^{\frac{-1-n}{n}}-c_n v_{k}^{\frac{-2}{n}}=\frac{c_n}{\omega_n ^{\frac{2}{n}}}a_{p_k }+\mathcal{O}_{p_k} ({v_k}^{\frac{1}{n}})$$
which completes the proof.
\end{Dem}

\begin{Cor} The solutions of the isoperimetric problem enclosing a small volume $v$ are pseudo-balls of constant mean curvature in a small neighborhood of the maxima of the scalar curvature function. 
\end{Cor} 

\begin{Res} \label{Cor1} Let $\mathcal{M}$ be a compact Riemannian manifold, let $$S=Sup_{\tilde{p}\in\mathcal{M}}\{ Sc(\tilde{p})\}.$$\\ 
                              Then\\
                              the isoperimetric profile $I_{\mathcal{M}}(v) $ has the following asymptotic expansion  
                              in a neighborhood of the origin:
                             \begin{equation}
I_{\mathcal{M}}(v)=c_n v^{\frac{n-1}{n}}\left( 1-\frac{S}{2n(n+2)}\left( \frac{v}{\omega_n}\right) ^{\frac{2}{n}} +o(v^{\frac{2}{n}})\right).
\end{equation}
\end{Res}

\begin{Dem}
By contradiction, applying lemma \ref{maxsc}.
\end{Dem}

Under stronger conditions, we can improve the remainder of this asymptotic expansion.

\begin{Res} \label{Cor2} Let $\mathcal{M}$ be a compact Riemannian manifold, let $$S:=Sup_{\tilde{p}\in\mathcal{M}}\{ Sc(\tilde{p})\}.$$ We assume that absolute maxima of $Sc$ are non degenerate critical points.\\ 
                              Then\\
                              the isoperimetric profile $I_{\mathcal{M}}(v) $ has the following asymptotic expansion in a neighborhood of the origin:
                             \begin{equation}
I_{\mathcal{M}}(v)=c_n v^{\frac{n-1}{n}}\left( 1-\frac{S}{2n(n+2)}\left( \frac{v}{\omega_n}\right) ^{\frac{2}{n}} +\mathcal{O}(v^{\frac{4}{n}})\right).
\end{equation}
\end{Res}

\begin{Dem}
            Let us set $\mathcal{N}_{r,p}^+ :=\{ exp_p (t\theta) | 0\leq t\leq r(1+x(p,r)(\theta))\} $.
            The critical points of $f(p,v):p\mapsto Vol_{n-1}(\beta(p,v))$ that are close to the maxima of the scalar curvature function, 
$p_1 , \cdots ,p_l$ are $C^{\infty}$ functions of $r$ denoted by $p_1 (r),\cdots ,p_l (r)$. This can be seen by an application of the implicit funcion theorem to the function $\nabla g$.\\
If $v$ is sufficiently small, $I_{\mathcal{M}}$ is achieved by a pseudo-ball $\beta(p_i (v),v)$ where
$\omega_n \rho(p,r)^n =v$ and the link between $r$ and $\rho$ is described at the end of section $1$. We set $f_i (v)=Vol_{n-1}(\beta(p_i (v),v))$, this is $C^{\infty}$ a function of $v^\frac{1}{n}$ and $$I_{\mathcal{M}}(v)=Min_{i\in \{ 1,\cdots , l\}}\left\lbrace f_i (v) \right\rbrace .$$\\ 
From lemma \ref{Ipexp} we deduce easily that
\begin{equation*}
I_{\mathcal{M}}(v)=c_n v^{\frac{n-1}{n}}\left( 1-\frac{1}{2n(n+2)}S\left( \frac{v}{\omega_n}\right) ^{\frac{2}{n}} +\mathcal{O}(v^{\frac{4}{n}})\right) .
\end{equation*}
\end{Dem}

Let $B$ be a ball in the model space of constant sectional curvature $K_0$.\\
 It is no difficult to check that for balls of small volume,
$$Vol(\partial B)=c_n Vol(B)^{\frac{n-1}{n}}\left[ 1 -\frac{1}{2n(n+2)}(n(n-1)K_0)\left( \frac{n}{\alpha_{n-1}}\right)^{\frac{2}{n}}Vol(B)^{\frac{2}{n}}\right] +\cdots $$
This permits to check that the expansion in theorem \ref{Cor1} coincides with theorem $1$ of \cite{Druet}.\\
To finish, we verify that, under the assumptions of the article \cite{Ye}, solutions of the isoperimetric problem for small volumes belong to the family of constant mean curvature spheres built by Ye. In order to make this possible we first show that this family coincides with the cmc pseudo-balls.\\
Here we use the notations of the paper \cite{Ye} and we use a trivialisation of the tangent bundle by an orthonormal frame field. 
\begin{Lemme}
           Let $p$ be a non degenerate critical point of the scalar curvature function of $\mathcal{M}$. Then 
           there exist $r_1$ so that for all $r<r_1$
           $$S_{r,\tau (r), r^2 \varphi (r), p}=\beta_1 (p(r),r)=\mathcal{N}_{r,p(r),x(r)}$$
           where $S_{r,\tau (r), r^2 \varphi (r), p}$ is constructed in the paper \cite{Ye} and is a parametrization of a foliation by cmc spheres that is constructed in the same article and $\beta_1 (p(r),r):=\beta(p(v(r)),v(r))$ where $ \beta$ is constructed in the preceding theorem.
\end{Lemme}

\begin{Dem} \\
1. The Riemannian center of mass $q(r)$ of the hypersurface $S_{r,\tau (r), r^2 \varphi (r), p}$ of Ye is close to $p$.

We apply lemma 1.4 to $p=p$ and $x(\theta)=r^{2}\phi(r)$. We find that $d(p,q(r))=r|A(r,p,r^{2}\phi(r))|$. This is due to the fact that $A(0,p,0)=0$, $|A(r,p,r^{2}\phi(r))|=O(r)$, from which we deduce $d(p,q(r))=O(r^{2})$.

2. The center of Ye $\exp_{p}(\tau(r))$ is close to $p$.
Indeed, following Ye, \cite{Ye} page 390, $d(p,\exp_{p}(\tau(r)))=O(r^{2})$.

3. The projection $\pi(r)$ of $S_{r,\tau (r), r^2 \varphi (r),p}$ on the sphere $S(q(r),r)$ is a  diffeomorphism $C^{\infty}$-close to the identity, when $r$ goes to $0$.

If $z$ is a point of $S_{r,\tau (r), r^2 \varphi (r),p}$, the angle between the geodesic segments from $z$ to  points $\exp_{p}(\tau(r))$ and $q(r)$ goes to $0$ uniformly in $r$ (we use a comparison theorem of   Riemannian geometry as in the proof of lemma 4.4 of \cite{Nar}).

4. We can write $S_{r,\tau (r), r^2 \varphi (r),p}$ as the graph of a function of the form $r(1+\tilde{x}(r))$ on the unit sphere on the tangent space in $q(r)$, and the $C^{2,\alpha}$ norm of $\tilde{x}(r)$ is small.

In fact, $r(1+\tilde{x}(r))$ is expressed as a function of $\phi(r)$ and of the inverse diffeomorphism of $\pi(r)$.

5. By construction, the function $\tilde{x}(r)$ belongs to the kernel of the linear differential operator $L$, and it satisfies the pseudo-ball equation, hence, applying lemma 1.6, $S_{r,\tau (r), r^2 \varphi (r), p}$ is the pseudo-ball $\beta_1 (q(r),r)$ of center of mass $q(r)$ and radius $r$.

6. $S_{r,\tau (r), r^2 \varphi (r),p }$ is a constant mean curvature surface, $q(r)$ is a critical point of the map $p\mapsto vol_{n-1}(\beta(p,v))$, where $v$ is the enclosed volume by $S_{r,\tau (r), r^2 \varphi (r),p}$. Then, $q(r)=p(r)$.\\
From definitions it is immediate to verify that  
$$S_{r,\tau (r), r^2 \varphi (r), p}=\beta_1 (p(r),v(r)).$$
\end{Dem}

\begin{Cor}
           If $\mathcal{M}$ is a compact manifold and the scalar curvature function of $\mathcal{M}$,   
           $Sc_{\mathcal{M}}$ has no degenerate maxima as critical points, then
           the solutions of the isoperimetric problem belong to the family constructed by Ye \cite{Ye} in a neighborhood of the non degenerate maxima of the scalar curvature function. 
\end{Cor}

\begin{Dem} Let $p_1 , \ldots , p_k$ be non degenerate maxima of the scalar curvature function, 
then, as the preceding lemma and theorem show, the solutions of the isoperimetric problem for small volume $v$ belong to the finite set $$\left\lbrace \beta_1(p_1(r),r),\ldots , \beta_1 (p_k (r),r)\right\rbrace$$ that coincides with $$\left\lbrace S_{r,\tau (r), r^2 \varphi (r),p_1},\ldots , S_{r,\tau (r), r^2 \varphi (r),p_k } \right\rbrace$$          
\end{Dem}

\textbf{Remark:} In the case of a compact manifold with non degenerate maxima of the scalar curvature we can continue the calculation of the asymptotic expansion to obtain the next non trivial coefficient.

      \newpage
      \markboth{References}{References}
      \bibliographystyle{alpha}
      \bibliography{these}
      \addcontentsline{toc}{section}{\numberline{}References}
      \emph{Stefano Nardulli\\ Dipartimento di Metodi e Modelli Matematici\\ Viale delle Scienze Edificio 8 - 90128 Palermo\\ email: nardulli@unipa.it} 
\end{document}